\documentclass[a4paper,english,oneside]{article}
\usepackage[T1]{fontenc}
\usepackage[utf8]{inputenc} 

\usepackage{geometry}
\usepackage{amsmath}
\usepackage{amsfonts}
\usepackage{amssymb}
\usepackage{amsthm}
\usepackage[english]{babel}
\usepackage{color}
\usepackage{lmodern}
\usepackage{pstricks-add}
\usepackage{dsfont}
\usepackage{hyperref}
\hypersetup{hidelinks}

\usepackage{mathscinet}
\usepackage{enumitem}
\usepackage[cal=boondoxo,bb=ams]{mathalfa}
\usepackage{microtype}

\usepackage{mathtools}
\usepackage{esint}

\numberwithin{equation}{section}

\theoremstyle{plain}
\newtheorem{theorem}{Theorem}[section]
\newtheorem{lemma}[theorem]{Lemma}

\theoremstyle{definition}
\newtheorem{definition}[theorem]{Definition}

\theoremstyle{remark}
\newtheorem{remark}{Remark}

    \DeclareMathOperator\supp{supp}

 \date{}
     
\begin{document}

\title{On a semilinear wave equation in anti-de Sitter spacetime: the critical case}

\author{Alessandro Palmieri$\,^{\mathrm{a}}$, Hiroyuki Takamura$\,^{\mathrm{b}}$}

\date{
\small{ $\,^\mathrm{a}$ Mathematical Institute, Tohoku University, Aoba, Sendai 980-8578, Japan} \\
\small{ $\,^\mathrm{b}$ Mathematical Institute/Research Alliance Center of Mathematical Sciences, Tohoku University, Aoba, Sendai 980-8578, Japan}  }

\maketitle
\begin{center}
\normalsize{\today}
\end{center}

\begin{abstract}
In the present paper we prove the blow-up in finite time for local solutions of a semilinear Cauchy problem associated with a wave equation in anti-de Sitter spacetime in the critical case. According to this purpose, we combine an ODI result with an iteration argument, by using an explicit integral representation formula for the solution to a linear Cauchy problem associated with the wave equation in anti-de Sitter spacetime in one space dimension.
\end{abstract}

\begin{flushleft}
\end{flushleft}

\begin{flushleft}
\textbf{AMS Classification (2020)} Primary: 35L05, 35L71, 35B44; Secondary: 33C05, 35A22, 35C15.

\end{flushleft}

\section{Introduction}

In the previous work \cite{PalTak22}, we studied the following semilinear Cauchy problem associated with a semilinear wave equation in anti-de Sitter spacetime
\begin{align}\label{anti deSitter semi} 
\begin{cases}
\partial_t^2 v- c^2\mathrm{e}^{2Ht} \Delta v+ b\partial_t v +m^2 v=f(t,v), & x\in \mathbb{R}^n, \ t\in (0,T), \\
v(0,x)= \varepsilon v_0(x), & x\in \mathbb{R}^n, \\
\partial_t v(0,x)= \varepsilon v_1(x), & x\in \mathbb{R}^n,
\end{cases}
\end{align} where $c,H$ are positive constants, $b,m^2$ are nonnegative real parameters satisfying $b^2\geqslant 4m^2$, $\varepsilon>0$ is a parameter describing the size of initial data,  $T=T(\varepsilon)\in (0,\infty]$ is the lifespan of a classical solution $v$ (i.e., the maximal existence time) and the nonlinear term is given by
\begin{align}\label{nonlinearity}
f(t,v)\doteq \Gamma(t) \left(\int_{\mathbb{R}^n}|v(t,y)|^p\mathrm{d}y\right)^\beta |v|^p,
\end{align} where $p>1$, $\beta\geqslant 0$, and $\Gamma=\Gamma(t)$ is a suitable positive function. 
In cosmology, the constant $H$ is the so-called Hubble constant, $m$ is the mass of a particle and $b$ is taken equal to the space dimension $n$ (cf. (0.6) in \cite{YagGal09}).

We recall that the assumption $b^2\geqslant 4m^2$ guarantees that the damping term $b\partial_t v$ is somehow dominant over the mass term $m^2 v$ (cf. \cite[Subsections 1.1-1.3]{ER18}).

More precisely, considering as $\Gamma$ factor  in \eqref{nonlinearity}
\begin{align}\label{def Gamma anti dS}
\Gamma(t)\doteq \mu \, \mathrm{e}^{\varrho t} (1+t)^{\varsigma},
\end{align} for some $\mu >0$ and $\varrho,\varsigma\in\mathbb{R}$,
 we have established the following threshold values for $\varrho$ 
 \begin{align}\label{def rho crit}
 \varrho_{\mathrm{crit}}(n,H,b,m^2,\beta,p)\doteq 
 \begin{cases}
 \vphantom{\Big(} \frac{1}{2}(b-\sqrt{b^2-4m^2})((\beta+1)p-1)+nH(\beta+1)(p-1) & \mbox{if} \ n\leqslant N, \\
 \vphantom{\Big(} \frac{1}{2}(b+nH)((\beta+1)p-1)+nH-(n-1)H(\beta+1)-\frac{H}{p} & \mbox{if} \ n> N,
 \end{cases}
 \end{align} where $$N=N(H,b,m^2,p)\doteq \tfrac{2}{p}+\tfrac{\sqrt{b^2-4m^2}}{H}.$$
 
 While for the case $n\leqslant N$ we provided a full picture on blow-up results for local weak solutions to \eqref{anti deSitter semi} for $\varrho\geqslant\varrho_{\mathrm{crit}}(n,H,b,m^2,\beta,p)$, cf. \cite[Theorems 1.6, 1.7 and 1.8]{PalTak22}, when $n >N$ only the case $\varrho>\varrho_{\mathrm{crit}}(n,H,b,m^2,\beta,p)$ was investigated (see \cite[Theorems 1.9]{PalTak22}).
 
 Aim of the present paper is to study the blow-up for local solutions to \eqref{anti deSitter semi} under suitable sign conditions for the Cauchy data in the threshold case $\varrho=\varrho_{\mathrm{crit}}(n,H,b,m^2,\beta,p)$ for $n> N$ and to derive the corresponding upper bound estimates for the lifespan.
 
 Furthermore, since we consider the limit case $\varrho=\varrho_{\mathrm{crit}}(n,H,b,m^2,\beta,p)$ we have to prescribe a lower bound for the power of the polynomial factor in \eqref{def Gamma anti dS}, namely, 
 \begin{align}
\varsigma_{\mathrm{crit}}(n,H,b,m^2,p)  &\doteq -\tfrac{1}{p}. \label{def sigma crit}
\end{align}

The threshold case that we treat in the present work is somehow a blow-up result for a critical case. Consequently, the approach that we use to prove the blow-up in finite time of the spatial average $V=V(t)$ of a local solution $v$ to \eqref{anti deSitter semi} is inspired by the one in \cite{TakWak11} for the derivation of the sharp upper bound estimate for the lifespan of a local solution to the semilinear wave equation in the critical case when $n\geqslant 4$. As in \cite{YZ06}, when $n\geqslant 2$ we use the Radon transform with respect to the spatial variables to handle the problem as it was in one space dimension.

The main difficulty in our argument will be the derivation of a sequence of lower bound estimates for the nonlinear term $\|v(t,\cdot)\|^p_{L^p(\mathbb{R}^n)}$. In particular, we will derive an iteration frame for $\|v(t,\cdot)\|^p_{L^p(\mathbb{R}^n)}$ combining two estimates involving the Radon transform of $v(t,\cdot)$. A fundamental tool for this kind of argument is provided by \emph{Yagdjian's integral transform approach}. Indeed, we will make use of an explicit integral representation formula for the solution to a linear one dimensional wave equation in anti-de Sitter spacetime in order to derive one of the aforementioned inequalities involving the Radon transform of $v(t,\cdot)$.

After deriving this sequence of lower bound estimates for $\|v(t,\cdot)\|^p_{L^p(\mathbb{R}^n)}$, we will derive in turn a sequence of lower bound estimates for $V$ with an additional polynomial factor. Combining these lower bound estimates for $V$ with a comparison argument for an ODE with ``critical'' exponential growth, we will be able to derive the desired upper bound estimates for the lifespan.

We point out that the speed of propagation, namely, the function $a(t)\doteq c\, \mathrm{e}^{Ht}$, is exponentially increasing in the previous semilinear wave equation. Moreover, the amplitude of the forward light-cone is given by $$A(t)\doteq \int_0^t a(\tau)\, \mathrm{d}\tau=\frac{c}{H}\left(\mathrm{e}^{Ht}-1\right).$$ 
This means that, considering smooth solutions, if we assume $v_0$ and $v_1$ compactly supported in $B_R\doteq \{x\in\mathbb{R}^n:|x|\leqslant R\}$, given a local solution $v$ to \eqref{anti deSitter semi}, we have that
\begin{align}\label{support condition sol}
\supp v(t,\cdot) \subset B_{R+A(t)} \ \ \mbox{for any} \ t\in (0,T).
\end{align} For the proof of this support condition one can use the property of finite speed of propagation or, alternatively, the explicit representation formulas from the series of works by Galstian and Yagdjian \cite{YagGal08,YagGal09a,YagGal09e}.

In this second part of the introduction, we provide a short summary of the results in the literature for wave models with a not-flat and time-dependent metric in the spacetime.

In the case of de Sitter spacetime, i.e. for $H<0$ in \eqref{anti deSitter semi}, the wave equation was considered by several authors. We recall the integral representation formulas (and their applications) established by  Yagdjian and Yagdjian-Galstian in \cite{YagGal09,Yag09,Yag10,Yag12,Yag13,Yag15,Yag19} and the global existence results for semilinear wave models in \cite{Na14} and \cite{ER18}.
Concerning blow-up results, we recall the blow-up result with a pure imaginary mass term in \eqref{anti deSitter semi}, namely, when we replace $m$ with $im$, both for de Sitter and anti-de Sitter spacetime in \cite[Proposition 1.1]{Na20}. Moreover, in \cite{TsuWa21} a blow-up result is proved in a de Sitter-type spacetime when $m=0$. Finally, in our recent paper \cite{PalTak22}, we slightly improved some result from \cite{Yag09} for the semilinear wave equation in de Sitter spacetime with the same nonlinear term as in \eqref{nonlinearity}, providing further the lifespan estimates (see \cite[Theorems 1.2,1.3 and 1.4]{PalTak22}).

For anti-de Sitter spacetime, we refer to  \cite{Gal03,YagGal08,YagGal09a,YagGal09e}, where among other things,  $L^p-L^q$ estimates are derived for the solutions to the corresponding linear Cauchy problem. 
Moreover, as we have already mentioned above, in \cite{PalTak22} some blow-up results for \eqref{anti deSitter semi} have been proved together with the corresponding upper bound estimates for the lifespan.

Finally, for the wave equation in Einstein-de Sitter spacetime (that is, for the d'Alember operator $\square_{\mathrm{EdS}}=\partial_t^2-t^{-2k}\Delta +b t^{-1}\partial_t$ with $k\in (0,1)$ and $b\geqslant 0$) we cite the papers \cite{GalKinYag10,GalYag14} for the linear model, \cite{GalYag20,TsuWa20,Pal20,Pal21,TsuWa21h,TsuWa22} for the semilinear model with power nonlinearity $|v|^p$ and \cite{HHP21,HHP21ar,TsuWa21G} for the semilinear model with nonlinearity of derivative type $|\partial_t v|^p$.

It is interesting to compare our approach in this paper to deal with a critical case in comparison to those in \cite{TsuWa20,Pal20,Pal21} for the treatment of the corresponding critical cases in Einsten-de Sitter spacetime. Indeed, while the critical case in \cite{TsuWa20} is studied by using the approach from \cite{ZH14}, in \cite{Pal20,Pal21} the method from \cite{WY19} is adapted to the case with time-dependent coefficients. However, in both cases the employment of the techniques from \cite{ZH14,WY19} to the case with time-dependent coefficients produces some restrictive assumptions: in \cite{TsuWa20} an upper bound is prescribed for the size of the ball containing the supports of the Cauchy data, whilst in in \cite{Pal20} a restriction on the multiplicative constant in the damping term is prescribed. We emphasize that with the approach of the present work, no restriction of this kind (either on the size of the support for the data or on the range for the multiplicative constants in the lower order terms) appears.


\subsection{Main results} \label{Subsection Main results}

Before stating the main theorem for \eqref{anti deSitter semi}, we recall the notion of weak solutions to \eqref{anti deSitter semi}  that has been employed in \cite{PalTak22}. We stress that, although we call them weak solutions, actually for these solutions more regularity than for usual distributional solutions is required with respect to the time variable, in order to handle a space average that is a $\mathcal{C}^2$ function with respect to $t$. Indeed, in \cite{PalTak22} we worked with the larger class of solutions that can be considered when employing the spatial average for proving the blow-up in finite time.

\begin{definition} \label{Def sol anti de Sitter} Let $v_0,v_1\in L^1_{\mathrm{loc}}(\mathbb{R}^n)$ such that $\supp v_0, \ \supp v_1 \subset B_R$ for some $R>0$. We say that $$v\in  \mathcal{C}^1\left([0,T), L^{1}_{\mathrm{loc}}(\mathbb{R}^n)\right) \ \mbox{such that} \  f(t,v) \in L^1_{\mathrm{loc}}((0,T)\times \mathbb{R}^n),$$ where the definition of the nonlinear term $f(t,u)$ is given in \eqref{nonlinearity},  is a \emph{weak solution} to \eqref{anti deSitter semi} on $[0,T)$ if  $v$ fulfills the support condition \eqref{support condition sol} and the integral identity
\begin{align}
& \int_{\mathbb{R}^n} \partial_t v(t,x) \varphi(t,x) \, \mathrm{d}x-\int_{\mathbb{R}^n} v(t,x) \varphi_t(t,x) \, \mathrm{d}x+ b\int_{\mathbb{R}^n}  v(t,x) \varphi(t,x) \, \mathrm{d}x \notag \\
& \qquad + \int_0^t\int_{\mathbb{R}^n} v(s,x) \left(\varphi_{ss}(s,x)-c^2\mathrm{e}^{2Hs}\Delta \varphi(s,x) - b \,\varphi_s(s,x) +m^2 \varphi(s,x) \right)  \mathrm{d}x \, \mathrm{d}s \notag \\
& \quad  = \varepsilon \int_{\mathbb{R}^n}  v_1(x) \varphi(0,x) \, \mathrm{d}x +\varepsilon \int_{\mathbb{R}^n} v_0(x) \big(b\, \varphi(0,x) -\varphi_t(0,x)\big) \, \mathrm{d}x \notag  \\ & \qquad + \int_0^t \Gamma(s) \left(\int_{\mathbb{R}^n}  |v(s,y)|^p  \, \mathrm{d}y\right)^\beta \int_{\mathbb{R}^n}  |v(s,x)|^p  \varphi(s,x) \, \mathrm{d}x \, \mathrm{d}s  \label{def weak sol int rel anti dS}
\end{align}
holds for any $t\in (0,T)$ and any test function $\varphi\in\mathcal{C}^\infty_0 ([0,T)\times \mathbb{R}^n)$.
\end{definition}

We point out explicitly that in the present paper we work with classical solutions to \eqref{anti deSitter semi} since we need to employ an integral representation formula that requires a pointwise evaluation of the Cauchy data and of the nonlinear term. Nonetheless, it is clear that classical solutions to \eqref{anti deSitter semi} are in particular weak solutions according to Definition \ref{Def sol anti de Sitter}. We emphasize that it was necessary to recall Definition \ref{Def sol anti de Sitter} since in what follows we are going to use some results from \cite{PalTak22} (see Section \ref{Section preliminary results}) that have been obtained for weak solutions in the aforementioned sense.

\begin{theorem} \label{Thm anti dS nlin lb poly growth}

Let $n\geqslant 1$ and $b,m^2\geqslant 0$ such that $b^2\geqslant 4m^2$. Let us assume $\beta\geqslant 0$ and  $p>1$ such that
\begin{align}\label{1st lb for V from V0 is dominant}
\frac{n}{2}-\frac{\sqrt{b^2-4m^2}}{2H}>\frac{1}{p}.
\end{align} 
 For $\varrho=\varrho_{\mathrm{crit}}(n,H,b,m^2,\beta,p)$ and $\varsigma> \varsigma_{\mathrm{crit}}(n,H,b,m^2,p)$, we consider 
\begin{align} \label{Gamma factor polynomial growth anti dS nlin lb}
\Gamma(t)\doteq \mu\,\mathrm{e}^{\varrho_{\mathrm{crit}}(n,H,b,m^2,\beta,p) t} (1+t)^\varsigma
\end{align} for some $\mu>0$ in the term $f(t,v)$ given by \eqref{nonlinearity}. 

 Let us assume that $(v_0,v_1)\in \mathcal{C}^2_0(\mathbb{R})\times  \mathcal{C}^1_0(\mathbb{R})$ 
  are nonnegative and nontrivial functions with supports contained in $ B_R$ for some $R>0$.
 Let $v\in \mathcal{C}^2\left([0,T)\times \mathbb{R}^n\right)$ be a classical solution to the Cauchy problem  \eqref{anti deSitter semi}  with lifespan $T=T(\varepsilon)$.

 Then, there exists a positive constant $\varepsilon_0=\varepsilon_0(n,c,H,b,m^2,\beta,p,\mu,\varsigma,v_0,v_1,R)$ such that for any $\varepsilon\in (0,\varepsilon_0]$ the classical solution $v$ blows up in finite time. Furthermore, the following upper bound estimates for the lifespan hold
\begin{align} \label{upper bound for the lifespan polynomial growth anti dS nlin}
T(\varepsilon) \leqslant  
\begin{cases}
 C \varepsilon^{-\frac{p((\beta+1)p-1)}{1+\varsigma p}} & \mbox{if} \ \ \varsigma \in \big(-\frac{1}{p},0\big], \\
 C \varepsilon^{-p((\beta+1)p-1)} & \mbox{if} \ \ \varsigma \in \big(0,\frac{1}{p}\big), \\
 C \varepsilon^{-\frac{(\beta+1)p-1}{\varsigma }} & \mbox{if} \ \ \varsigma \in \big[\frac{1}{p},\infty\big),
\end{cases}
\end{align} where the positive constant $C$ is independent of $\varepsilon$.
\end{theorem}

The next sections of the paper are organized as follows: in Section \ref{Section preliminary results} we recall briefly some estimates derived in \cite{PalTak22}; in Section 3 we prove Theorem \ref{Thm anti dS nlin lb poly growth}. More precisely, this proof is split in the following intermediate steps: in Subsection \ref{Subsection ODI} we derive a comparison argument for an \emph{ordinary differential inequality} (ODI) with ``critical'' exponential growth; in Subsection \ref{Subsubsection Yagdjian repr formula} we derive an integral representation formula for the solution to the one dimensional linear Cauchy problem associated with the wave equation in anti-de Sitter spacetime; in Subsection \ref{Subsection iteration frame ||v(t)||^p L^p} we derive the crucial iteration frame for $\|v(t,\cdot)\|^p_{L^p(\mathbb{R}^n)}$ and in Subsection \ref{Subsection sequence lb Lp norm v(t)} we use this iteration frame to derive a sequence of lower bound estimates for $\|v(t,\cdot)\|^p_{L^p(\mathbb{R}^n)}$; hence, we use in turn such a sequence to derive a sequence of lower bound estimates for the spatial average of the local solution in Subsection \ref{Subsection sequence lb V(t)} and complete the proof in the case $\varsigma\in (\varsigma_{\mathrm{crit}}(n,H,b,m^2,p),0]$; finally, in Subsection \ref{Subsubsection varsigma >0} we conclude the proof also for the case $\varsigma>0$.

\section{Preliminary results} \label{Section preliminary results}

Before beginning with the proof of Theorem \ref{Thm anti dS nlin lb poly growth}, we recall some estimates that are proved in Section 3 of \cite{PalTak22}.

\subsection{Iteration frame for the spatial average}\label{Section iteration frame anti-dS}

Let $v$ be a local classical solution to \eqref{anti deSitter semi}. In particular, by using the property of finite speed of propagation, we have that $v$ satisfies the support condition \eqref{support condition sol}. We set
\begin{align*}
V(t)\doteq \int_{\mathbb{R}^n} v(t,x) \, \mathrm{d}x \qquad \mbox{for} \ t\in(0,T).
\end{align*}
In \cite[Subsection 3.1]{PalTak22}, 
 we proved the identity
\begin{align}\label{V''+bV'+m2 V}
V''(t)+bV'(t)+m^2 V(t)= \Gamma(t) \left(\int_{\mathbb{R}^n}|v(t,x)|^p\mathrm{d}x\right)^{\beta+1}.
\end{align}
We underline that \eqref{V''+bV'+m2 V} is obtained by choosing a suitable cut-off function in \eqref{def weak sol int rel anti dS}, that localizes the forward light-cone on the strip $[0,t]\times\mathbb{R}^n$.
Hence, by factorizing the differential operator $\partial_t^2+b\partial_t +m^2$ in \eqref{V''+bV'+m2 V}, we derived then the following inequality for $V$:
\begin{align}\label{fundamental ineq V(t)}
V(t)\geqslant  \mathrm{e}^{-\alpha_2 t} \int_0^t \mathrm{e}^{(\alpha_2-\alpha_1) s} \int_0^s \mathrm{e}^{\alpha_1 \tau} \, \Gamma(\tau)\left(\|v(\tau,\cdot)\|^{p}_{L^p(\mathbb{R}^n)} \right)^{(\beta+1)} \, \mathrm{d}\tau  \, \mathrm{d}s,
\end{align} for $t\geqslant 0$, where $\alpha_1, \alpha_2$ are the roots of the quadratic equation $\alpha^2-b\alpha+m^2=0$.

The estimate in \eqref{fundamental ineq V(t)} is very important since it will be used in Subsection \ref{Subsection sequence lb V(t)} to derive a sequence of lower bound estimates for $V(t)$ from the sequence of lower bound estimates for $\|v(t,\cdot)\|^p_{L^p(\mathbb{R}^n)}$  derived in Subsection \ref{Subsection sequence lb Lp norm v(t)}.

\subsection{First lower bound estimate for $\|v(t,\cdot)\|^p_{L^p(\mathbb{R}^n)}$}

In Subsection 3.2 of \cite{PalTak22}, by working with a weighted spatial average of  $v$, with the weight function given by a suitable positive solution of the linear adjoint equation with separable variables, we derived the following lower bound estimates 
\begin{align} \label{||v(t)||^p Lp lower bound step -1}
\|v(t,\cdot)\|^p_{L^p(\mathbb{R}^n)} & \geqslant \widetilde{B} \, \varepsilon^p \mathrm{e}^{\left[-\frac{1}{2}(b+H)p+(n-1)H(1-\frac{p}{2})\right]t}
\end{align} for $t\geqslant 0$ and for a suitable constant $\widetilde{B}=\widetilde{B}(n,c,H,b,m^2,p,v_0,v_1)>0$.

\section{Proof of Theorem \ref{Thm anti dS nlin lb poly growth}}

\label{Subsection anti dS nlin lb poly growth}

The proof of the results in the critical case $\varrho=\varrho_{\mathrm{crit}}(n,H,b,m^2,\beta,p)$ when we are in the case $n>N$ is more delicate than the ones for $n\leqslant N$ seen in \cite{PalTak22}.  Roughly speaking, the functional $V$ is no longer sufficient to show the blow-up in finite time of a local solution. The tools that we are going to use in this section are inspired by the ones for treatment of the critical case for the semilinear wave equation in the flat case. In particular, we are going to combine the approaches from \cite{YZ06} and \cite{TakWak11} with some ideas from \cite{PalRei19} for the treatment of a semilinear wave equation with time-dependent coefficients.

\subsection{ODI comparison argument in the critical case} \label{Subsection ODI}

We state and prove a Kato-type lemma for exponentially growing functions in the ``critical case''. This result is the counterpart in our framework of Lemma 2.1 in \cite{TakWak11}.

\begin{lemma} \label{Lemma ODI} Let $b,m^2$ be nonnegative real numbers such that $b^2\geqslant 4m^2$. We set $$\alpha_1\doteq \tfrac b2+ \tfrac{1}{2}\sqrt{b^2-4m^2} \quad \mbox{and} \quad \alpha_2\doteq \tfrac b2- \tfrac{1}{2}\sqrt{b^2-4m^2}.$$ Let $q>1$, $k_0,k_1\in\mathbb{R}$ such that 
\begin{align}
& k_0+(q-1)k_1=0 \label{critical condition k0, k1}, \\
&  k_1+ \alpha_1 \geqslant 0 \label{condition k0, k1}.
\end{align}
Suppose that $G\in\mathcal{C}^2([0,T))$ satisfies
\begin{align}
& G''(t)+b\, G'(t)+m^2G(t)\geqslant B \, \mathrm{e}^{k_0 t}|G(t)|^q \qquad \mbox{for} \ t\geqslant 0, \label{ODI for G statement} \\
& G(t)\geqslant K \mathrm{e}^{k_1 t} \qquad \qquad \qquad \qquad \qquad \qquad \quad \ \, 	\, \mbox{for} \ t\geqslant T_0,  \label{lower bound G statement} \\
& G(0),G'(0)\geqslant 0 \quad \mbox{and} \quad \alpha_1 G(0)+G'(0)>0 \label{Initial condition G statement}
\end{align} with $T_0\in[0,T)$ and for some  positive constants $B,K$.
Let us define 
\begin{equation} \label{def T1 and K0}
\begin{split}
T_1 & \doteq 
\begin{cases}
\max\left\{T_0,\tilde{T}_0, (k_1+\alpha_1)^{-1}\right\} & \mbox{if} \ k_1+\alpha_1>0, \\
\max\left\{T_0,\tilde{T}_0\right\}  & \mbox{if} \  k_1+\alpha_1=0,
\end{cases} \\
 K_0 & \doteq 
 \begin{cases} 
 \displaystyle{\left(\frac{q+1}{B}\right)^{\frac{1}{q-1}} \left(\frac{k_1+\alpha_1}{1-\mathrm{e}^{-\vartheta}}\right)^{\frac{2}{q-1}}} & \mbox{if} \ k_1+\alpha_1>0, \\
\displaystyle{\left(\frac{q+1}{B}\right)^{\frac{1}{q-1}} \left(\kappa\vartheta \right)^{-\frac{2}{q-1}}} & \mbox{if} \ k_1+\alpha_1=0,
 \end{cases}
\end{split}
\end{equation} 
where $\vartheta\in (0,\frac{q-1}{2})$  and $\kappa\in (0,T_0)$ are arbitrarily chosen and 
\begin{align*}
\widetilde{T}_0= \widetilde{T}_0 (b,m^2,G(0),G'(0))\doteq \begin{cases} \frac{1}{\alpha_1-\alpha_2}\ln \left(1+\frac{(\alpha_1-\alpha_2)G(0)}{\alpha_1 G(0)+G'(0)}\right) & \mbox{if} \ \  b^2>4m^2,\\
\frac{G(0)}{\tfrac{b}{2} G(0)+G'(0)}    & \mbox{if} \ \ b^2=4m^2.\end{cases}
\end{align*}  If the multiplicative constant on the right-hand side of \eqref{lower bound G statement} satisfies $K\geqslant K_0$, then, the lifespan of $G$ is finite and fulfills $T\leqslant  2T_1$.
\end{lemma}

\begin{remark}
 Since $q>1$, the two conditions \eqref{critical condition k0, k1} and \eqref{condition k0, k1} imply immediately that
\begin{align}
 k_0-\alpha_1(q-1)  \leqslant  0 \label{condition k0, k1 n2}.
\end{align} We will see that also the condition in \eqref{condition k0, k1 n2} on $k_0,k_1$, along with \eqref{critical condition k0, k1} and \eqref{condition k0, k1}, is fundamental for the proof of Lemma \ref{Lemma ODI}.
\end{remark}

\begin{remark}\label{Remark k1 >= alpha 2}
As we are going to see in the proof of Lemma \ref{Lemma ODI}, we may assume without loss of generality that the coefficient $k_1$ appearing in the lower bound for $G$ in \eqref{lower bound G statement} satisfies $$k_1\geqslant -\alpha_2 \geqslant - \alpha_1,$$ where the equality $-\alpha_2= -\alpha_1$ holds only for $b^2=4m^2$. Indeed, it is possible to replace the lower bound for $G$ in \eqref{lower bound G statement} with \eqref{G >= G lin} below in order to get $k_1\geqslant -\alpha_2$.

In particular, it makes sense to consider the limit case $k_1+\alpha_1=0$ (and, consequently, to modify accordingly $T_1$ and $K_0$ as in \eqref{def T1 and K0}) only in the balanced case $b^2=4m^2$.

We point out explicitly, that the condition \eqref{condition k0, k1} in our application of Lemma \ref{Lemma ODI} will be always satisfied thanks to \eqref{1st lb for V from V0 is dominant}.
\end{remark}

\begin{proof}
 By contradiction, we assume that $G(t)$ is defined for any $t\in [0,2T_1]$. We will show that this fact is not compatible with the choice $K\geqslant K_0$.

Let us begin by proving that $G$ is actually nonnegative for any $t\in[0,T).$ By using the factorization 
\begin{align}\label{factorization G''+bG'+m2G proof}
\mathrm{e}^{-\alpha_2 t} \frac{\mathrm{d}}{\mathrm{d} t} \left( \mathrm{e}^{(\alpha_2-\alpha_1) t} \frac{\mathrm{d}}{\mathrm{d} t} \left( \mathrm{e}^{\alpha_1 t}G(t)\right) \right) =G''(t)+b\, G'(t)+m^2G(t)\geqslant 0, 
\end{align}
straightforward computations  lead to the following lower bound estimate for $G(t)$:
\begin{align}\label{G >= G lin}
G(t)\geqslant G_{\mathrm{lin}}(t)\doteq 
\begin{cases}
\vphantom{\Big(} \displaystyle{ \frac{\alpha_1\, \mathrm{e}^{-\alpha_2 t}-\alpha_2\,\mathrm{e}^{-\alpha_1 t}}{\alpha_1-\alpha_2} G(0)+  \,\frac{\mathrm{e}^{-\alpha_2 t}-\mathrm{e}^{-\alpha_1 t}}{\alpha_1-\alpha_2}G'(0)} & \mbox{if} \ \alpha_1 \neq \alpha_2,  \\
\vphantom{\bigg(} \displaystyle{\left(1+\tfrac{b}{2}t\right)\mathrm{e}^{-\frac{b}{2} t} G(0) + t\,  \mathrm{e}^{-\frac{b}{2} t} G'(0)} & \mbox{if} \ \alpha_1=\alpha_2.  
\end{cases}
\end{align}
See \cite[Section 2.1]{PalTak22} for a detailed derivation of the previous kind of inequality for a function $G$ satisfying the ordinary differential inequality in \eqref{ODI for G statement}.

Let us introduce now the further time-dependent function $F(t)\doteq \mathrm{e}^{\alpha_1 t} G(t)$. By the previous considerations it results $F(t)\geqslant \mathrm{e}^{\alpha_1 t} G_{\mathrm{lin}}(t)\geqslant 0$. From \eqref{ODI for G statement} it follows that 
\begin{align}\label{ODI for F proof}
F''(t)+(b-2\alpha_1)F'(t)\geqslant B \mathrm{e}^{(k_0-\alpha_1(q-1))t} (F(t))^q.
\end{align} By using \eqref{factorization G''+bG'+m2G proof}, which can be rewritten as 
\begin{align*}
\mathrm{e}^{-\alpha_2 t} \frac{\mathrm{d}}{\mathrm{d} t} \left( \mathrm{e}^{(\alpha_2-\alpha_1) t} F'(t) \right)\geqslant 0, 
\end{align*} we get easily that 
\begin{align} \label{lower bound F'}
F'(t)\geqslant \mathrm{e}^{(\alpha_1-\alpha_2) t} F'(0) = \mathrm{e}^{(\alpha_1-\alpha_2) t} (G'(0)+\alpha_1 G(0))>0
\end{align} thanks to the last sign assumption for the initial values of $G$ in \eqref{Initial condition G statement}. Therefore, we multiply both sides of the inequality in \eqref{ODI for F proof} by $F'(t)$, arriving at 
\begin{align}
\frac{1}{2}\frac{\mathrm{d}}{\mathrm{d} t} \left( \left(F'(t)\right)^2 \right) &\geqslant (\alpha_1-\alpha_2)\left(F'(t)\right)^2 + \frac{B \mathrm{e}^{(k_0-\alpha_1(q-1))t} }{q+1} \frac{\mathrm{d}}{\mathrm{d} t} \left((F(t))^{q+1}\right) \notag \\
&\geqslant  \frac{B \mathrm{e}^{(k_0-\alpha_1(q-1))t} }{q+1} \frac{\mathrm{d}}{\mathrm{d} t} \left((F(t))^{q+1}\right), \label{neglected quadratic term}
\end{align} where in the second estimate we used the fact that $\alpha_1\geqslant \alpha_2$. Integrating both sides of the previous inequality over $[0,t]$, we obtain
\begin{align*}
 \left(F'(t)\right)^2- \left(F'(0)\right)^2 &\geqslant \frac{2B }{q+1}  \int_0^t \mathrm{e}^{(k_0-\alpha_1(q-1))\tau} \frac{\mathrm{d}}{\mathrm{d} \tau} \left((F(\tau))^{q+1}\right) \mathrm{d}\tau \\
  &\geqslant \frac{2B \mathrm{e}^{(k_0-\alpha_1(q-1))t}}{q+1}  \left[(F(t))^{q+1}-(F(0))^{q+1}\right], 
\end{align*} where in the second step we used \eqref{condition k0, k1 n2}. By straightforward computations, we find
\begin{align}
F(t)-F(0)\geqslant \mathrm{e}^{\alpha_1 t}G_{\mathrm{lin}}(t)-G(0) = \begin{cases}
\vphantom{\Big(} \displaystyle{ (\alpha_1G(0)+G'(0))\frac{\mathrm{e}^{(\alpha_1-\alpha_2) t}-1}{\alpha_1-\alpha_2} } & \mbox{if} \ \alpha_1 \neq \alpha_2,  \\
\vphantom{\bigg(} \displaystyle{ \left(\tfrac{b}{2}G(0)+G'(0)\right) t} & \mbox{if} \ \alpha_1=\alpha_2,  
\end{cases} \label{lower bound G>G lin proof}
\end{align}  so, in both cases it holds $F(t)>F(0)\geqslant 0$ for $t>0$.
Hence,
\begin{align}
 \left(F'(t)\right)^2 & \geqslant \frac{2B \mathrm{e}^{(k_0-\alpha_1(q-1))t}}{q+1} (F(t))^{q} \left[F(t)-F(0)\left(\frac{F(0)}{F(t)}\right)^{q}\right]\notag \\
 & \geqslant \frac{2B \mathrm{e}^{(k_0-\alpha_1(q-1))t}}{q+1} (F(t))^{q} \left[F(t)-F(0)\right]. \label{lower bound (F')^2}
\end{align} Next we show that $F(t)\geqslant 2F(0)$ for $t\geqslant \widetilde{T}_0$. For $\alpha_1\neq \alpha_2$, from \eqref{lower bound G>G lin proof} we have
\begin{align*}
F(t)-2F(0) \geqslant  (\alpha_1 G(0)+G'(0)) \frac{\mathrm{e}^{(\alpha_1-\alpha_2) t}-1}{\alpha_1-\alpha_2} -G(0),
\end{align*} and then $F(t)\geqslant 2F(0)$ is satisfied provided that
\begin{align*}
\frac{\mathrm{e}^{(\alpha_1-\alpha_2) t}-1}{\alpha_1-\alpha_2}  \geqslant \frac{G(0)}{\alpha_1 G(0)+G'(0)} \qquad \Longleftrightarrow \qquad t\geqslant  \frac{1}{\alpha_1-\alpha_2}\ln \left(1+\frac{(\alpha_1-\alpha_2)G(0)}{\alpha_1 G(0)+G'(0)}\right) = \widetilde{T}_0.
\end{align*}  For $\alpha_1= \alpha_2$, the situation is even simpler since 
\begin{align*}
F(t)-2F(0)\geqslant \left(\tfrac{b}{2}G(0)+G'(0)\right) t -G(0).
\end{align*} Hence, for $t\geqslant \widetilde{T}_0$ it follows from \eqref{lower bound (F')^2} that
\begin{align*}
\left(F'(t)\right)^2 & \geqslant \frac{2B \mathrm{e}^{(k_0-\alpha_1(q-1))t}}{q+1} (F(t))^{q} \left[F(t)-F(0)\right] \\
 & \geqslant \frac{B \mathrm{e}^{(k_0-\alpha_1(q-1))t}}{q+1} (F(t))^{q+1}. 
\end{align*} From the last inequality we get
\begin{align*}
F'(t) & \geqslant \sqrt{\frac{B }{q+1}} \ \mathrm{e}^{\frac{1}{2}(k_0-\alpha_1(q-1))t}(F(t))^{\frac{q+1}{2}}.
\end{align*} Now we multiply both sides of the previous inequality by $(F(t))^{-(1+\vartheta)}$ for some $\vartheta\in (0,\tfrac{q-1}{2})$, obtaining
\begin{align*}
\frac{F'(t)}{(F(t))^{1+\vartheta}}= \frac{\mathrm{d}}{\mathrm{d}t}\left(-\frac{1}{\vartheta (F(t))^\vartheta}\right) \geqslant \sqrt{\frac{B }{q+1}} \ \mathrm{e}^{\frac{1}{2}(k_0-\alpha_1(q-1))t}(F(t))^{\frac{q-1}{2}-\vartheta}.
\end{align*} Integrating both sides of the last inequality over $[T_1,t]$, it results
\begin{align}\label{fundamental estimate F proof}
\frac{1}{\vartheta}\left(\frac{1}{(F(T_1))^{\vartheta}}-\frac{1}{(F(t))^{\vartheta}}\right) \geqslant   \sqrt{\frac{B }{q+1}} \int_{T_1}^t \mathrm{e}^{\frac{1}{2}(k_0-\alpha_1(q-1))\tau}(F(\tau))^{\frac{q-1}{2}-\vartheta} \mathrm{d}\tau.
\end{align} The next step will consist in determining a suitable upper (resp. lower) bound estimate for the left-hand (resp. right-hand) side of \eqref{fundamental estimate F proof}. In order to derive both these estimates, we employ \eqref{lower bound G statement} to establish the following lower bound estimate for $F$
\begin{align} \label{lower bound F proof}
F(t)\geqslant K \mathrm{e}^{(\alpha_1+k_1)t} \qquad \mbox{for} \ \ t\geqslant T_0.
\end{align} Therefore, employing \eqref{lower bound F proof}, for $t\geqslant T_1$ we find, on the one hand,
\begin{align} \label{fundamental estimate F proof n2}
\frac{1}{\vartheta}\left(\frac{1}{(F(T_1))^{\vartheta}}-\frac{1}{(F(t))^{\vartheta}}\right) < \frac{1}{\vartheta}\frac{1}{(F(T_1))^{\vartheta}} \leqslant \frac{\mathrm{e}^{-\vartheta(k_1+\alpha_1)T_1}}{\vartheta K^\vartheta}.
\end{align}
In this last part of the proof, we have to consider separately the case $k_1+\alpha_1>0$ and the case $k_1+\alpha_1=0$. For $k_1+\alpha_1>0$ we have
\begin{align}
 \sqrt{\frac{B }{q+1}} &\int_{T_1}^t \mathrm{e}^{\frac{1}{2}(k_0-\alpha_1(q-1))\tau}(F(\tau))^{\frac{q-1}{2}-\vartheta} \mathrm{d}\tau  \notag \\
& \quad \geqslant \sqrt{\frac{B }{q+1}} \frac{K^{\frac{q-1}{2}}}{K^{\vartheta}}\int_{T_1}^t \exp\left\{\left[\tfrac{1}{2}(k_0+(q-1)k_1)-\vartheta(k_1+\alpha_1)\right]\tau\right\} \mathrm{d}\tau \notag \\
& \quad = \sqrt{\frac{B }{q+1}} \frac{K^{\frac{q-1}{2}}}{K^{\vartheta}}\int_{T_1}^t \exp\left(-\vartheta(k_1+\alpha_1)\tau\right) \mathrm{d}\tau \notag \\
& \quad = \sqrt{\frac{B }{q+1}} \frac{K^{\frac{q-1}{2}}}{\vartheta K^{\vartheta}}(k_1+\alpha_1)^{-1} \left(\mathrm{e}^{-\vartheta(k_1+\alpha_1)T_1}-\mathrm{e}^{-\vartheta(k_1+\alpha_1)t}\right),
\label{fundamental estimate F proof n3}
\end{align} where in the second step we used the threshold condition \eqref{critical condition k0, k1}. So, combining \eqref{fundamental estimate F proof}, \eqref{fundamental estimate F proof n2} and \eqref{fundamental estimate F proof n3} and using the fact that $T_1\geqslant (k_1+\alpha_1)^{-1}$, we have for $t=2T_1$
\begin{align*}
1 &>\sqrt{\frac{B }{q+1}} K^{\frac{q-1}{2}} (k_1+\alpha_1)^{-1} \left(1-\mathrm{e}^{-\vartheta(k_1+\alpha_1)(t-T_1)}\right)
\\ &\geqslant \sqrt{\frac{B }{q+1}} K^{\frac{q-1}{2}} (k_1+\alpha_1)^{-1} \left(1-\mathrm{e}^{-\vartheta}\right) =\left(\frac{K}{K_0}\right)^{\frac{q-1}{2}},
\end{align*} which contradicts the condition $K\geqslant K_0$  for $k_1+\alpha_1>0$. 
On the other hand, for $k_1+\alpha_1=0$ we have
\begin{align}
 \sqrt{\frac{B }{q+1}} \int_{T_1}^t \mathrm{e}^{\frac{1}{2}(k_0-\alpha_1(q-1))\tau}(F(\tau))^{\frac{q-1}{2}-\vartheta} \mathrm{d}\tau 
& \quad \geqslant \sqrt{\frac{B }{q+1}} \frac{K^{\frac{q-1}{2}}}{K^{\vartheta}}\int_{T_1}^t \exp\left(-\vartheta(k_1+\alpha_1)\tau\right) \mathrm{d}\tau \notag \\
& \quad = \sqrt{\frac{B }{q+1}} \frac{K^{\frac{q-1}{2}}}{ K^{\vartheta}}(t-T_1),
\label{fundamental estimate F proof n4}
\end{align} Analogously as before, we combine \eqref{fundamental estimate F proof}, \eqref{fundamental estimate F proof n2} and \eqref{fundamental estimate F proof n4}, obtaining for $t=\kappa+T_1\leqslant 2T_1$ with $\kappa\in (0,T_0)$
\begin{align*}
1 &>\sqrt{\frac{B }{q+1}} K^{\frac{q-1}{2}} \kappa \vartheta
 =\left(\frac{K}{K_0}\right)^{\frac{q-1}{2}},
\end{align*} which contradicts the condition $K\geqslant K_0$ for $k_1+\alpha_1=0$.
This completes the proof.
\end{proof}

\begin{remark} In the previous proof, we neglected in \eqref{neglected quadratic term} the influence of the term $(\alpha_1-\alpha_2) (F'(t))^2$ in the intermediate steps that lead to \eqref{lower bound (F')^2}. If we used \eqref{lower bound F'} to estimate this term from below and we kept the resulting term till the estimate \eqref{lower bound (F')^2}, by having ignored the other term which appears on the right-hand side of \eqref{lower bound (F')^2}, we would obtain $(F'(t))^2 \gtrsim \mathrm{e}^{2(\alpha_1-\alpha_2)t}$ that would imply in turn $F(t)\gtrsim  \mathrm{e}^{(\alpha_1-\alpha_2)t}$ for $t$ sufficiently large. In particular, by comparing this lower bound for $F$ with the one in \eqref{lower bound F proof}, we see that the latter is stronger provided that $k_1\geqslant -\alpha_2$.

As we pointed out in Remark \ref{Remark k1 >= alpha 2}, we may always assume without loss of generality that $k_1\geqslant -\alpha_2$ and, consequently, that  \eqref{lower bound F proof} provides the best lower bound estimate for $F$.
\end{remark}

\subsection{Integral representation formula for the 1-d linear Cauchy problem} \label{Subsubsection Yagdjian repr formula}

We derive now an integral representation formula for the solution of the following linear inhomogeneous Cauchy problem in one space dimension
\begin{align}\label{anti deSitter lin 1d} 
\begin{cases}
\partial_t^2 v- c^2\mathrm{e}^{2Ht} \partial_x^2 v+ b\partial_t v +m^2 v=g(t,x), & x\in \mathbb{R}, \ t>0, \\
v(0,x)= v_0(x), & x\in \mathbb{R}, \\
\partial_t v(0,x)= v_1(x), & x\in \mathbb{R}.
\end{cases}
\end{align}
Similar representation formulas are well-established in the literature by Yagdjian and Galstian-Yagdjian (cf. the introduction) both in de Sitter and anti-de Sitter spacetime for the Klein-Gordon equation, normalizing $c=1=H$. By means of a suitable change of variables first and, then, employing the so-called dissipative transformation, we will easily derive a representation formula for the solution of \eqref{anti deSitter lin 1d} by using \cite[Theorems 3-4]{YagGal09a}.

Let us consider the change of variables
\begin{align} \label{change of variables from t,x to tau,y}
\tau \doteq Ht, \quad y\doteq \tfrac Hc x.
\end{align}
Elementary computations show that $v$ solves the following Cauchy problem with respect to $(\tau,y)$ 
\begin{align}\label{anti deSitter lin 1d tau,y} 
\begin{cases}
\partial_\tau^2 v- \mathrm{e}^{2\tau} \partial_y^2 v+ \frac{b}{H}\partial_\tau v +\frac{m^2}{H^2} v=\frac{1}{H^2}g(\tfrac{\tau}{H},\tfrac{cy}{H}), & y\in \mathbb{R}, \ \tau>0, \\
v(0,y)= v_0(\tfrac{cy}{H}), & y\in \mathbb{R}, \\
\partial_\tau v(0,x)= \frac{1}{H} v_1(\tfrac{cy}{H}), & y\in \mathbb{R}.
\end{cases}
\end{align}
Applying the transformation $v(\tau,y)\doteq \mathrm{e}^{-\frac{b\tau }{2H}}w(\tau,y)$, we find that $v$ solves the Cauchy problem in \eqref{anti deSitter lin 1d tau,y} if and only if $w$ solves 
\begin{align}\label{anti deSitter lin 1d w} 
\begin{cases}
\partial_\tau^2 w- \mathrm{e}^{2\tau} \partial_y^2 w-\frac{b^2-4m^2}{4H^2} w=\frac{\mathrm{e}^{\frac{b \tau}{2H}}}{H^2}g(\tfrac{\tau}{H},\tfrac{cy}{H}), & y\in \mathbb{R}, \ \tau>0, \\
w(0,y)= v_0(\tfrac{cy}{H}), & y\in \mathbb{R}, \\
\partial_\tau w(0,x)= \frac{b}{2H}v_0(\tfrac{cy}{H})+\frac{1}{H} v_1(\tfrac{cy}{H}), & y\in \mathbb{R}.
\end{cases}
\end{align}

Let us set $G(\tau,y)\doteq \frac{\mathrm{e}^{\frac{b\tau}{2H}}}{H^2}g(\tfrac{\tau}{H},\tfrac{cy}{H})$, $w_0(y)\doteq v_0(\tfrac{cy}{H})$, $w_1(y)\doteq \frac{b}{2H}v_0(\tfrac{cy}{H})+\frac{1}{H} v_1(\tfrac{cy}{H})$ and $\nu\doteq \frac{\sqrt{b^2-4m^2}}{2H}$. In particular, with the terminology from \cite{YagGal09}, $w$ satisfies a Klein-Gordon equation in anti-de Sitter spacetime with complex-valued curved mass $\nu$.

According to Theorems 3 and 4 from \cite{YagGal09a}, $w$ admits the following representation
\begin{align}
w(\tau,y) & = \int_{0}^\tau \int_{y-\mathrm{e}^\tau+\mathrm{e}^\sigma}^{y+\mathrm{e}^\tau-\mathrm{e}^\sigma} \widetilde{E}(\tau,y;\sigma,y_0;\nu) G(\sigma,y_0)\, \mathrm{d}y_0 \,\mathrm{d}\sigma \notag \\
& \quad +\tfrac{1}{2} \, \mathrm{e}^{-\frac{\tau}{2}} \big(w_0(y+\mathrm{e}^\tau-1)+w_0(y-\mathrm{e}^\tau+1)\big) \notag \\
& \quad +  \int_{y-\mathrm{e}^\tau+1}^{y+\mathrm{e}^\tau-1} \widetilde{K}_0(\tau,y;y_0;\nu) w_0(y_0)\, \mathrm{d}y_0 + \int_{y-\mathrm{e}^\tau+1}^{y+\mathrm{e}^\tau-1} \widetilde{K}_1(\tau,y;y_0;\nu) w_1(y_0)\, \mathrm{d}y_0, \label{representation for w}
\end{align} with the kernel functions given by
\begin{align}
 \widetilde{E}(\tau,y;\sigma,y_0;\nu) & \doteq 2^{-2\nu} \mathrm{e}^{-\nu(\tau+\sigma)}\left((\mathrm{e}^{\tau}+\mathrm{e}^\sigma)^2-(y-y_0)^2)\right)^{-\frac{1}{2}+\nu} \notag \\
 & \qquad \times\mathrm{F}\left(\frac{1}{2}-\nu,\frac{1}{2}-\nu;1;\frac{(\mathrm{e}^{\tau}-\mathrm{e}^\sigma)^2-(y-y_0)^2}{(\mathrm{e}^{\tau}+\mathrm{e}^\sigma)^2-(y-y_0)^2}\right), \label{def E tilde} \\
  \widetilde{K}_0(\tau,y;y_0;\nu) & \doteq  -\frac{\partial}{\partial \sigma}\widetilde{E}(\tau,y;\sigma,y_0;\nu) \Big|_{\sigma=0}, \label{def K0 tilde}\\
  \widetilde{K}_1(\tau,y;y_0;\nu) & \doteq  \widetilde{E}(\tau,y;0,y_0;\nu), \label{def K1 tilde}
\end{align} 
 where $\mathrm{F}(\frac{1}{2}-\nu,\frac{1}{2}-\nu;1;\cdot)$ denotes the Gauss hypergeometric function.

Inverting the change of variables in \eqref{change of variables from t,x to tau,y}, we may rewrite the four addends in the representation for $w$ in \eqref{representation for w} in a more convenient way. Let us begin with the double integral involving the source term
\begin{align}
& \int_{0}^\tau \int_{y-\mathrm{e}^\tau+\mathrm{e}^\sigma}^{y+\mathrm{e}^\tau-\mathrm{e}^\sigma} \widetilde{E}(\tau,y;\sigma,y_0;\nu) G(\sigma,y_0)\, \mathrm{d}y_0 \,\mathrm{d}\sigma \notag \\
& \qquad = \frac{1}{H^2}\int_{0}^\tau \int_{y-\mathrm{e}^\tau+\mathrm{e}^\sigma}^{y+\mathrm{e}^\tau-\mathrm{e}^\sigma} \widetilde{E}(\tau,y;\sigma,y_0;\nu) \, \mathrm{e}^{\frac{b\sigma}{2H}}g(\tfrac{\sigma}{H},\tfrac{c y_0}{H})\, \mathrm{d}y_0 \,\mathrm{d}\sigma  \notag\\
& \qquad = \frac{1}{c}\int_{0}^{\tfrac{\tau}{H}} \int_{\tfrac{c}{H}(y-\mathrm{e}^\tau+\mathrm{e}^{Hs})}^{\tfrac{c}{H}(y+\mathrm{e}^\tau-\mathrm{e}^{Hs})} \widetilde{E}(\tau,y;Hs,\tfrac{Hz}{c};\nu) \, \mathrm{e}^{\frac{b}{2}s}g(s,z)\, \mathrm{d}z \,\mathrm{d}s \notag\\
& \qquad = \frac{1}{c}\int_{0}^{t} \int_{x-\tfrac{c}{H}(\mathrm{e}^{Ht}-\mathrm{e}^{Hs})}^{x+\tfrac{c}{H}(\mathrm{e}^{Ht}-\mathrm{e}^{Hs})} \mathrm{e}^{\frac{b}{2}s} \widetilde{E}(Ht,\tfrac{Hx}{c};Hs,\tfrac{Hz}{c};\nu) \, g(s,z)\, \mathrm{d}z \,\mathrm{d}s. \label{w double intgral}
\end{align} By \eqref{change of variables from t,x to tau,y} and the definition for $w_0$, we can easily express the second addend in \eqref{representation for w} as follows
\begin{align}
\tfrac{1}{2}\mathrm{e}^{-\frac{\tau}{2}} \big(w_0(y+\mathrm{e}^\tau-1)+w_0(y-\mathrm{e}^\tau+1)\big) &= \tfrac{1}{2}\mathrm{e}^{-\frac{H}{2}t} \big(v_0(x+\tfrac{c}{H}(\mathrm{e}^{Ht}-1))+v_0(x-\tfrac{c}{H}(\mathrm{e}^{Ht}-1))\big) \notag\\ &= \tfrac{1}{2}\mathrm{e}^{-\frac{H}{2}t} \big(v_0(x+A(t))+v_0(x-A(t))\big). \label{w travelling wave}
\end{align} Finally, we rewrite together the remaining integrals in \eqref{representation for w}. By using the definition of $w_0,w_1$ and \eqref{change of variables from t,x to tau,y}, we get
\begin{align*}
   &\int_{y-\mathrm{e}^\tau+1}^{y+\mathrm{e}^\tau-1} \left(\widetilde{K}_0(\tau,y;y_0;\nu) w_0(y_0)+ \widetilde{K}_1(\tau,y;y_0;\nu) w_1(y_0)\right) \mathrm{d}y_0 \\
   & \quad= \int_{y-\mathrm{e}^\tau+1}^{y+\mathrm{e}^\tau-1} \left(\widetilde{K}_0(\tau,y;y_0;\nu) v_0(\tfrac{cy_0}{H}) + \tfrac{1}{H}\widetilde{K}_1(\tau,y;y_0;\nu) \left(\tfrac{b}{2}v_0(\tfrac{cy_0}{H})+ v_1(\tfrac{cy_0}{H})\right)\right) \mathrm{d}y_0 \\
   & \quad = \frac{1}{c}\int_{\tfrac{c}{H}(y-\mathrm{e}^\tau+1)}^{\tfrac{c}{H}(y+\mathrm{e}^\tau-1)} \left(\left(H\widetilde{K}_0(\tau,y;\tfrac{Hz}{c};\nu)+ \tfrac{b}{2} \widetilde{K}_1(\tau,y;\tfrac{Hz}{c};\nu)\right) v_0(z) + \widetilde{K}_1(\tau,y;\tfrac{Hz}{c};\nu)  v_1(z)\right) \mathrm{d}z \\
   & \quad = \frac{1}{c}\int_{x-\tfrac{c}{H}(\mathrm{e}^{Ht}-1)}^{x+\tfrac{c}{H}(\mathrm{e}^{Ht}-1)} \left(H\widetilde{K}_0(Ht,\tfrac{Hx}{c};\tfrac{Hz}{c};\nu)+ \tfrac{b}{2} \widetilde{K}_1(Ht,\tfrac{Hx}{c};\tfrac{Hz}{c};\nu)\right) v_0(z) \mathrm{d}z  \\& \quad \qquad + \frac{1}{c}\int_{x-\tfrac{c}{H}(\mathrm{e}^{Ht}-1)}^{x+\tfrac{c}{H}(\mathrm{e}^{Ht}-1)} \widetilde{K}_1(Ht,\tfrac{Hx}{c};\tfrac{Hz}{c};\nu)  v_1(z) \mathrm{d}z.
\end{align*} We analyze more carefully the kernel for the first data, that is
\begin{align*}
 & H\widetilde{K}_0(Ht,\tfrac{Hx}{c};\tfrac{Hz}{c};\nu)+ \tfrac{b}{2} \widetilde{K}_1(Ht,\tfrac{Hx}{c};\tfrac{Hz}{c};\nu) \\
 & \quad =  -H \frac{\partial}{\partial \sigma}\widetilde{E}(Ht,\tfrac{Hx}{c};\sigma,\tfrac{Hz}{c};\nu) \Big|_{\sigma=0} +\tfrac{b}{2} \widetilde{E}(Ht,\tfrac{Hx}{c};0,\tfrac{Hz}{c};\nu) \\
 & \quad =  -\frac{\partial}{\partial s}\widetilde{E}(Ht,\tfrac{Hx}{c};Hs,\tfrac{Hz}{c};\nu) \Big|_{s=0} +\tfrac{b}{2} \widetilde{E}(Ht,\tfrac{Hx}{c};0,\tfrac{Hz}{c};\nu) \\
 & \quad =  -\frac{\partial}{\partial s} \left(\mathrm{e}^{\frac{b}{2}s}\widetilde{E}(Ht,\tfrac{Hx}{c};Hs,\tfrac{Hz}{c};\nu) \right)\Big|_{s=0} +b \widetilde{E}(Ht,\tfrac{Hx}{c};0,\tfrac{Hz}{c};\nu),
\end{align*} where in the second equality we used the chain rule $H \frac{\partial}{\partial \sigma}\big|_{\sigma=0}= \frac{\partial}{\partial s}\big|_{s=0}$ for $\sigma=Hs$. Consequently,
\begin{align}
   &\int_{y-\mathrm{e}^\tau+1}^{y+\mathrm{e}^\tau-1} \left(\widetilde{K}_0(\tau,y;y_0;\nu) w_0(y_0)+ \widetilde{K}_1(\tau,y;y_0;\nu) w_1(y_0)\right) \mathrm{d}y_0 \notag \\
   & \ \ = \frac{1}{c}\int_{x-\tfrac{c}{H}(\mathrm{e}^{Ht}-1)}^{x+\tfrac{c}{H}(\mathrm{e}^{Ht}-1)} \left(-\partial_s \left(\mathrm{e}^{\frac{b}{2}s}\widetilde{E}(Ht,\tfrac{Hx}{c};Hs,\tfrac{Hz}{c};\nu) \right)\Big|_{s=0} +b \widetilde{E}(Ht,\tfrac{Hx}{c};0,\tfrac{Hz}{c};\nu)\right) v_0(z) \mathrm{d}z  \notag \\& \quad \qquad + \frac{1}{c}\int_{x-\tfrac{c}{H}(\mathrm{e}^{Ht}-1)}^{x+\tfrac{c}{H}(\mathrm{e}^{Ht}-1)} \widetilde{E}(Ht,\tfrac{Hx}{c};0,\tfrac{Hz}{c};\nu) v_1(z) \mathrm{d}z. \label{w single integrals}
\end{align}
Using the inverse transformation $v(t,x)=\mathrm{e}^{-\frac{b}{2}t}w(\tau,y)$ in \eqref{representation for w} and combining \eqref{w double intgral}, \eqref{w travelling wave} and \eqref{w single integrals}, we conclude
\begin{align*}
v(t,x)& = \frac{1}{c}\int_{0}^{t} \int_{x-\tfrac{c}{H}(\mathrm{e}^{Ht}-\mathrm{e}^{Hs})}^{x+\tfrac{c}{H}(\mathrm{e}^{Ht}-\mathrm{e}^{Hs})} \mathrm{e}^{-\frac{b}{2}(t-s)} \widetilde{E}(Ht,\tfrac{Hx}{c};Hs,\tfrac{Hz}{c};\nu) \, g(s,z)\, \mathrm{d}z \,\mathrm{d}s 
\\& \quad +\frac{1}{2}\,\mathrm{e}^{-\frac{b+H}{2}t} \big[v_0\big(x+\tfrac{c}{H}(\mathrm{e}^{Ht}-1)\big)+v_0\big(x-\tfrac{c}{H}(\mathrm{e}^{Ht}-1)\big)\big] \\
&\quad +\frac{1}{c}\int_{x-\tfrac{c}{H}(\mathrm{e}^{Ht}-1)}^{x+\tfrac{c}{H}(\mathrm{e}^{Ht}-1)} \left(-\tfrac{\partial}{\partial s}+b\right) \left(\mathrm{e}^{-\frac{b}{2}(t-s)}\widetilde{E}(Ht,\tfrac{Hx}{c};Hs,\tfrac{Hz}{c};\nu) \right)\Big|_{s=0}  v_0(z) \mathrm{d}z  \notag \\ & \quad  + \frac{1}{c}\int_{x-\tfrac{c}{H}(\mathrm{e}^{Ht}-1)}^{x+\tfrac{c}{H}(\mathrm{e}^{Ht}-1)} \mathrm{e}^{-\frac{b}{2}t} \widetilde{E}(Ht,\tfrac{Hx}{c};0,\tfrac{Hz}{c};\nu) v_1(z) \mathrm{d}z.
\end{align*}

Summarizing the computations that we have just done and using the definitions \eqref{def E tilde}, \eqref{def K0 tilde} and \eqref{def K1 tilde}, we proved the following result.

\begin{lemma} Let $n=1$ and $c,H>0$, $b,m^2\geqslant 0$ such that $b^2\geqslant 4m^2$. Let us assume $v_0\in\mathcal{C}^2(\mathbb{R})$, $v_1\in\mathcal{C}^1(\mathbb{R})$ and $g\in\mathcal{C}([0,\infty),\mathcal{C}^1(\mathbb{R}))$. Then, the classical solution $v$ to the Cauchy problem \eqref{anti deSitter lin 1d} can be represented as follows
\begin{align*}
v(t,x)& =\int_{0}^{t} \int_{x-A(t)+A(s)}^{x+A(t)-A(s)} E(t,x;s,z;c,H,b,m^2) \, g(s,z)\, \mathrm{d}z \,\mathrm{d}s 
\\& \quad +\frac{1}{2}\,\mathrm{e}^{-\tfrac{(b+H)t}{2}} \big[v_0(x+A(t))+v_0(x-A(t))\big] \\
&\quad +\int_{x-A(t)}^{x+A(t)} K_0(t,x;z;c,H,b,m^2) v_0(z) \, \mathrm{d}z   + \int_{x-A(t)}^{x+A(t)} K_1(t,x;z;c,H,b,m^2) v_1(z) \,\mathrm{d}z,
\end{align*} where $A(t)=\frac{c}{H}(\mathrm{e}^{Ht}-1)$, the kernel functions are given by 
\begin{align}
E(t,x;s,z;c,H,b,m^2) & \doteq \frac{1}{H} \left(\frac{2c}{H}\right)^{-2\nu} \mathrm{e}^{-(\frac{b}{2}+\nu H)t}  	\, \mathrm{e}^{(\frac{b}{2}-\nu H)s}\left(\left(\tfrac{c}{H}\left(\mathrm{e}^{Ht}+\mathrm{e}^{Hs}\right)\right)^2-(x-z)^2\right)^{-\frac{1}{2}+\nu} \notag \\
 & \qquad \times\mathrm{F}\left(\frac{1}{2}-\nu,\frac{1}{2}-\nu;1;\frac{\left(\tfrac{c}{H}(\mathrm{e}^{Ht}-\mathrm{e}^{Hs})\right)^2-(x-z)^2}{\left(\tfrac{c}{H}(\mathrm{e}^{Ht}+\mathrm{e}^{Hs})\right)^2-(x-z)^2}\right), \label{def E} \\
  K_0(t,x;z;c,H,b,m^2) & \doteq  -\frac{\partial}{\partial s}E(t,x;s,z;c,H,b,m^2)\Big|_{s=0}+bE(t,x;0,z;c,H,b,m^2), \label{def K0}\\
  K_1(t,x;z;c,H,b,m^2) & \doteq  E(t,x;0,z;c,H,b,m^2), \label{def K1}
\end{align} and the parameter in the hypergeometric function is $\nu=\frac{\sqrt{b^2-4m^2}}{2H}$.
\end{lemma}

\begin{remark} In \eqref{def E} we could replace $\nu$ with $-\nu$, thanks to the following property of the hypergeometric function
\begin{align*}
\mathrm{F}(a_1,a_2,b;\zeta)=(1-\zeta)^{b-(a_1+a_2)} \mathrm{F}(b-a_1,b-a_2,b;\zeta)
\end{align*} see, for example, \cite[Equation (15.8.1)]{OLBC10}. However, we preferred the definition provided in \eqref{def E} since in this way we have no singular behavior for the hypergeometric function as $\zeta\to 1^+$ when $\nu>0$.
\end{remark}

\begin{remark}\label{Remark sign K0}
It is clear that the kernel functions $E$ and $K_1$ are nonnegative on the forward light-cone and on the base of the forward light-cone, respectively.   Now we want to show  that $K_0 (t,x;z;c,H,b,m^2)\geqslant 0$ for any $z$ such that $|x-z|\leqslant \tfrac cH (\mathrm{e}^{Ht}-1)$. 

Let us express in a more explicit way the kernel function $K_0$ in \eqref{def K0}. For the sake of brevity, we introduce the notation
\begin{align*}
\mathcal{E} (t,x;s,z;c,H,b,m^2)\doteq \mathrm{e}^{(\frac{b}{2}-\nu H)s}\left(\left(\tfrac{c}{H}\left(\mathrm{e}^{Ht}+\mathrm{e}^{Hs}\right)\right)^2-(x-z)^2\right)^{-\frac{1}{2}+\nu} \mathrm{F}\left(\tfrac{1}{2}-\nu,\tfrac{1}{2}-\nu;1;\zeta\right),
\end{align*} where 
\begin{align*}
\zeta=\zeta(t,x;s,z;c,H)\doteq \frac{\left(\tfrac{c}{H}(\mathrm{e}^{Ht}-\mathrm{e}^{Hs})\right)^2-(x-z)^2}{\left(\tfrac{c}{H}(\mathrm{e}^{Ht}+\mathrm{e}^{Hs})\right)^2-(x-z)^2}.
\end{align*} Clearly,
\begin{align*}
 K_0(t,x;z;c,H,b,m^2) & = \frac{1}{H} \left(\frac{2c}{H}\right)^{-2\nu} \mathrm{e}^{-(\frac{b}{2}+\nu H)t}  \left( -\tfrac{\partial}{\partial s}\big|_{s=0}+b\right)\mathcal{E}(t,x;0,z;c,H,b,m^2).
\end{align*}
Direct computations show that
\begin{align*}
\frac{\partial \mathcal{E}}{\partial s} & (t,x;s,z;c,H,b,m^2)  = (\tfrac b2 -\nu H) \mathcal{E} (t,x;s,z;c,H,b,m^2) \\ & +\frac{(2\nu-1)H (\tfrac cH)^2 (\mathrm{e}^{Ht}+\mathrm{e}^{Hs})\mathrm{e}^{Hs}}{\left(\tfrac{c}{H}(\mathrm{e}^{Ht}+\mathrm{e}^{Hs})\right)^2-(x-z)^2}\mathcal{E} (t,x;s,z;c,H,b,m^2) \\
&  +(\tfrac{1}{2}-\nu)^2 \mathrm{e}^{(\frac{b}{2}-\nu H)s}\left(\left(\tfrac{c}{H}\left(\mathrm{e}^{Ht}+\mathrm{e}^{Hs}\right)\right)^2-(x-z)^2\right)^{-\frac{1}{2}+\nu} \mathrm{F}\left(\tfrac{3}{2}-\nu,\tfrac{3}{2}-\nu;2;\zeta\right) \frac{\partial \zeta}{\partial s},
\end{align*} where we used the recursive identity $\mathrm{F}'(a_1,a_2;b;\zeta)=\frac{a_1 a_2}{b} \mathrm{F}(a_1+1,a_2+1;b+1;\zeta)$ for the derivative of the hypergeometric function (cf. \cite[Equation (15.5.1)]{OLBC10}).

Since
\begin{align*}
\frac{\partial \zeta}{\partial s}(t,x;s,z;c,H) &= 4 \tfrac{c^2}{H} \mathrm{e}^{H(t+s)}\frac{(x-z)^2+\left(\tfrac{c}{H}\right)^2 (-\mathrm{e}^{2Ht}+\mathrm{e}^{2Hs})}{\left[\left(\tfrac{c}{H}(\mathrm{e}^{Ht}+\mathrm{e}^{Hs})\right)^2-(x-z)^2\right]^2} \\
& \leqslant 4 \tfrac{c^2}{H} \mathrm{e}^{H(t+s)}\frac{2\left(\tfrac{c}{H}\right)^2 \mathrm{e}^{Hs} (\mathrm{e}^{Hs}-\mathrm{e}^{Ht})}{\left[\left(\tfrac{c}{H}(\mathrm{e}^{Ht}+\mathrm{e}^{Hs})\right)^2-(x-z)^2\right]^2}  \leqslant0
\end{align*} for $|x-z|\leqslant \tfrac cH (\mathrm{e}^{Ht}-\mathrm{e}^{Hs})$ and $s\in [0,t]$, we find
\begin{align}
 & -\frac{\partial \mathcal{E}}{\partial s}  (t,x;s,z;c,H,b,m^2) +b\mathcal{E} (t,x;s,z;c,H,b,m^2) \notag \\
& \quad \geqslant (\tfrac b2 +\nu H) \mathcal{E} (t,x;s,z;c,H,b,m^2)  +\frac{(1-2\nu)H (\tfrac cH)^2 (\mathrm{e}^{Ht}+\mathrm{e}^{Hs})\mathrm{e}^{Hs}}{\left(\tfrac{c}{H}(\mathrm{e}^{Ht}+\mathrm{e}^{Hs})\right)^2-(x-z)^2}\mathcal{E} (t,x;s,z;c,H,b,m^2).\label{mathcal E low bound}
\end{align}
For $\nu \in [0,\frac{1}{2}]$ the right-hand side of \eqref{mathcal E low bound} is nonnegative, so $K_0 (t,x;z;c,H,b,m^2)\geqslant 0$ for $|x-z|\leqslant \tfrac cH (\mathrm{e}^{Ht}-1)$. On the other hand, for $\nu>\tfrac 12$ we use in \eqref{mathcal E low bound} the upper bound estimate
\begin{align*}
\frac{(\tfrac cH)^2 (\mathrm{e}^{Ht}+\mathrm{e}^{Hs})\mathrm{e}^{Hs}}{\left(\tfrac{c}{H}(\mathrm{e}^{Ht}+\mathrm{e}^{Hs})\right)^2-(x-z)^2}\leqslant \frac{ (\mathrm{e}^{Ht}+\mathrm{e}^{Hs})\mathrm{e}^{Hs}}{4\mathrm{e}^{H(t+s)}} \leqslant \frac{ (1+\mathrm{e}^{-H(t-s)})}{4}\leqslant \frac{1}{2}
\end{align*} for $|x-z|\leqslant \tfrac cH (\mathrm{e}^{Ht}-\mathrm{e}^{Hs})$ and $s\in [0,t]$, to derive the following lower bound estimate
\begin{align*}
  -\frac{\partial \mathcal{E}}{\partial s}  (t,x;s,z;c,H,b,m^2) & +b\mathcal{E} (t,x;s,z;c,H,b,m^2)  \\
& \quad \geqslant (\tfrac b2 +\nu H+ \tfrac 12 (1-2\nu)H) \mathcal{E} (t,x;s,z;c,H,b,m^2) \\
& \quad \geqslant \tfrac 12 (b+H) \mathcal{E} (t,x;s,z;c,H,b,m^2)\geqslant 0
\end{align*} for $|x-z|\leqslant \tfrac cH (\mathrm{e}^{Ht}-\mathrm{e}^{Hs})$ and $s\in [0,t]$. For $s=0$ the previous inequality provides the nonnegativity of $K_0$ on the domain of integration.
\end{remark}

\subsection{Iteration frame for $\|v(t,\cdot)\|^p_{L^p(\mathbb{R}^n)}$ via the Radon transform} \label{Subsection iteration frame ||v(t)||^p L^p}

The idea to apply the Radon transform to reduce somehow the problem to a one-dimensional one when $n\geqslant 2$ was introduced in the study of the critical case for the semilinear wave equation in the flat case in \cite{YZ06}. Here we will follow the main ideas from \cite[Section 5]{PalRei19} to derive an iteration frame for the nonlinear term $\|v(t,\cdot)\|^p_{L^p(\mathbb{R}^n)}$. In particular, the representation formula obtained via Yagdjian's integral transform approach will have a crucial role in the explicit representation of the Radon transform of $v$. We emphasize that the case $n=1$ can be considered as well; however, rather than working with the Radon transform of $v$, it is sufficient to work simply with $v$ (cf. Remark \ref{Remark n=1} below for further details).

We begin with the following remark: without loss of generality we may assume that a local solution is radially symmetric with respect to $x$. Indeed, when $v$ is not radial it possible to consider instead $$\bar{v}(t,r)\doteq \fint_{\mathbb{S}^{n-1}}v(t,r\omega) \mathrm{d}\sigma_{\omega} \qquad \mbox{for} \ t\in [0,T),	\ r\geqslant 0. $$ Let us clarify the meaning of this statement. By Jensen's inequality we have
\begin{align*}
\fint_{\mathbb{S}^{n-1}}|v(t,r\omega)|^p \mathrm{d}\sigma_{\omega}\geqslant \left|\fint_{\mathbb{S}^{n-1}}v(t,r\omega) \mathrm{d}\sigma_{\omega}\right|^p,
\end{align*} that is, $\overline{|v|^p}\geqslant |\bar{v}|^p$.
 Similarly, by Fubini's theorem we have
\begin{align*}
\int_{\mathbb{R}^n}|v(t,y)|^p\mathrm{d}y = \omega_n \int_0^\infty \fint_{\mathbb{S}^{n-1}}|v(t,r\omega)|^p \mathrm{d}\sigma_{\omega} r^{n-1} \mathrm{d}r \geqslant  \omega_n \int_0^\infty |\bar{v}(t,r)|^p r^{n-1} \mathrm{d}r = \|\bar{v}(t,\cdot)\|^p_{L^p(\mathbb{R}^n)},
\end{align*} where $\omega_n$ denotes the $(n-1)$-dimensional measure of the unit sphere of $\mathbb{R}^n$.

Consequently, combining the previous inequality, from \eqref{nonlinearity} we have
\begin{align*}
\overline{f(t,v)} = \Gamma(t) \left(\int_{\mathbb{R}^n}|v(t,y)|^p\mathrm{d}y\right)^\beta \overline{|v|^p}\geqslant \Gamma(t) \|\bar{v}(t,\cdot)\|^{\beta p}_{L^p(\mathbb{R}^n)} |\bar{v}|^p=f(t,\bar{v}).
\end{align*} Since the fundamental solution $E$ defined in \eqref{def E} is nonnegative on the forward light-cone and the averages with respect to the space variables of $v$ and $\bar{v}$ are equal,  the inequality $\overline{f(t,v)}\geqslant f(t,\bar{v})$, that we just proved,  allows us to assume without loss of generality that $v$ is radially symmetric in the proof of the blow-up result. 

Let us recall the definition of Radon transform of $v(t,\cdot)$ when $n\geqslant 2$.  Given $\rho\in\mathbb{R}$ and $\xi\in \mathbb{R}^n$, $|\xi|=1$, we have
\begin{align}\label{Radon transform}
\mathcal{R}[v](t,\rho,\xi)\doteq \int_{\{x\in\mathbb{R}^n: \, x\cdot \xi=\rho\}}v(t,x) \mathrm{d}\sigma_x= \int_{\{x\in\mathbb{R}^n:\, x\cdot \xi=0\}}v(t,\rho\xi+x) \mathrm{d}\sigma_x,
\end{align} where $\mathrm{d}\sigma_x$ is the Lebesgue measure on the corresponding hyperplanes. Since $v(t,\cdot)$ is radially symmetric with respect to $x$, it turns out that $\mathcal{R}[v]$ does not depend actually on $\xi$ and, moreover,
\begin{align}\label{Radon transform radial case}
\mathcal{R}[v](t,\rho)= \omega_{n-1}\int_{|\rho|}^\infty v(t,r)(r^2-\varrho^2)^{\frac{n-3}{2}} r\mathrm{d}r.
\end{align} 
Indeed, using polar coordinates $x=\omega r_1$ with $r_1=|x|$ and $\omega\in\mathbb{S}^{n-1}$ such that $\omega\cdot \xi=0$, we obtain
\begin{align*}
\mathcal{R}[v](t,\rho,\xi) & = \int_{\{x\in\mathbb{R}^n:\, x\cdot \xi=0\}}v(t,\rho\xi+x) \mathrm{d}\sigma_x = \int_0^\infty \int_{\substack{\omega\in\mathbb{S}^{n-1} \\ \omega\cdot \xi=0 \hphantom{a \,}}} v\left(t,\sqrt{\smash[b]{\rho^2+r_1^2}}\right) \mathrm{d}\sigma_\omega r_1^{n-2} \, \mathrm{d}r_1 \\
& = \omega_{n-1} \int_0^\infty  v\left(t,\sqrt{\smash[b]{\rho^2+r_1^2}}\right)  r_1^{n-2} \,  \mathrm{d}r_1.
\end{align*} Hence, employing the change of variables $r=(\rho^2+r_1^2)^{1/2}$ in the last integral, we find \eqref{Radon transform radial case}.

Using the identity $\mathcal{R}[\Delta v]=\partial_\rho^2 \mathcal{R}[v]$ and the linearity of $\mathcal{R}$, we find that $\mathcal{R}[v]$ satisfies the following one-dimensional Cauchy problem
\begin{align*}
\begin{cases}
\partial_t^2 \mathcal{R}[v]- c^2\mathrm{e}^{2Ht} \partial_\rho^2 \mathcal{R}[v]+ b\partial_t \mathcal{R}[v] +m^2 \mathcal{R}[v]=\mathcal{R}[f(t,v)], & \rho\in \mathbb{R}, \ t\in (0,T), \\
\mathcal{R}[v](0,\rho)= \varepsilon\mathcal{R}[v_0](\rho), & \rho\in \mathbb{R}, \\
\partial_t \mathcal{R}[v](0,\rho)= \varepsilon \mathcal{R}[v_1](\rho), & \rho\in \mathbb{R}.
\end{cases}
\end{align*}
From Subsection \ref{Subsubsection Yagdjian repr formula}, it follows the following integral representation formula 
\begin{align*}
\mathcal{R}[v](t,\rho) =\varepsilon (\mathcal{R}[v])_{\mathrm{lin}}(t,\rho)+  (\mathcal{R}[v])_{\mathrm{nlin}}(t,\rho),
\end{align*} with
\begin{align*}
 (\mathcal{R}[v])_{\mathrm{lin}}(t,\rho) & \doteq 
\frac{1}{2}\,\mathrm{e}^{-\tfrac{(b+H)t}{2}} \left(\mathcal{R}[v_0](\rho+A(t))+\mathcal{R}[v_0](\rho-A(t))\right) \\
&\quad +\int_{\rho-A(t)}^{\rho+A(t)} K_0(t,\rho;\eta;c,H,b,m^2) \mathcal{R}[v_0](\eta)\mathrm{d}\eta   \\
&\quad + \int_{\rho-A(t)}^{\rho+A(t)} K_1(t,\rho;\eta;c,H,b,m^2) \mathcal{R}[v_1](\eta) \mathrm{d}\eta,\\
  (\mathcal{R}[v])_{\mathrm{nlin}}(t,\rho)  & \doteq \int_{0}^{t} \int_{\rho-A(t)+A(s)}^{\rho+A(t)-A(s)} E(t,\rho;s,\eta;c,H,b,m^2) \, \mathcal{R}[f(\cdot,v)](s,\eta)\, \mathrm{d}\eta \,\mathrm{d}s,
\end{align*} where the definitions of the kernel functions are given in \eqref{def E}, \eqref{def K0} and \eqref{def K1}, respectively. 
Thanks to the assumptions on the Cauchy data in the statement of Theorem \ref{Thm anti dS nlin lb poly growth} 
and Remark \ref{Remark sign K0}, we have that $(\mathcal{R}[v])_{\mathrm{lin}}$ is a nonnegative function. Therefore, it results
\begin{align*}
\mathcal{R}[v](t,\rho)\geqslant (\mathcal{R}[v])_{\mathrm{nlin}}(t,\rho) .
\end{align*} We point out that $\mathcal{R}$ acts only on the factors in the nonlinear term $f(t,v)$ that depend on the space variable, that is,
\begin{align*}
\mathcal{R}[f(\cdot,v)](t,\rho)= \Gamma(t) \|v(t,\cdot)\|^{\beta p}_{L^p(\mathbb{R}^n)}  \mathcal{R}[|v|^p](t,\rho).
\end{align*} Thus,
\begin{align}\label{first lower bound Rv}
\mathcal{R}[v](t,\rho)\geqslant \int_{0}^{t} \Gamma(s) \|v(s,\cdot)\|^{\beta p}_{L^p(\mathbb{R}^n)} \int_{\rho-A(t)+A(s)}^{\rho+A(t)-A(s)} E(t,\rho;s,\eta;c,H,b,m^2) \,  \mathcal{R}[|v|^p](s,\eta)\, \mathrm{d}\eta \,\mathrm{d}s.
\end{align} From the support condition \eqref{support condition sol} for $v$, it follows that $$\mathrm{supp} 	\, \mathcal{R}[v](t,\cdot)\subset [-(R+A(t)),R+A(t)] \quad \mbox{for any} \ t\in [0,T).$$ Indeed, for $|\rho|>R+A(t)$ from the second representation in \eqref{Radon transform} we have that
\begin{align*}
\mathcal{R}[v](t,\rho,\xi) = \int_{\{x\in\mathbb{R}^n:\, x\cdot \xi=0\}}v(t,\rho\xi+x) \mathrm{d}\sigma_x =0,
\end{align*} due to the fact that on the hyperplane where we are integrating it holds
\begin{align*}
|\rho \xi+x|^2=\rho^2+|x|^2\geqslant \rho^2 \quad \Longrightarrow \quad |\rho \xi+x|> R+ A(t)
\end{align*} and, consequently, the considered hyperplane as empty intersection with the support of $v$.

In a completely analogous way, we have that $\mathrm{supp} 	\, \mathcal{R}[|v|^p](t,\cdot)\subset [-(R+A(t)),R+A(t)]$ for any $t\in [0,T).$

In the next step, we shrink the domain of integration with respect to $s$ in \eqref{first lower bound Rv} so that the support of $ \mathcal{R}[|v|^p](s,\eta)$ is a subset of the $\eta$-domain of integration. In other words, we look for $s\in [0,t]$ such that 
\begin{align*}
 [-(R+A(s)),R+A(s)] \subset [\rho - & A(t)+A(s) ,\rho + A(t)-A(s)] \\ &\quad \Longleftrightarrow  \quad 2A(s)\leqslant A(t)-|\rho|-R \\
 & \quad \Longleftrightarrow  \quad  s\leqslant  s_0(t,\rho,R)\doteq A^{-1}\left(\tfrac{1}{2}(A(t)-|\rho|-R)\right).
\end{align*} We point out that $s_0\geqslant 0$ if and only if $|\rho |\leqslant A(t) -R$.

Therefore, for $|\rho |\leqslant A(t) -R$ we obtain from \eqref{first lower bound Rv}
\begin{align} \label{second lower bound Rv}
\mathcal{R}[v](t,\rho)\geqslant \int_{0}^{s_0} \Gamma(s) \|v(s,\cdot)\|^{\beta p}_{L^p(\mathbb{R}^n)} \int_{-(R+A(s))}^{R+A(s)} E(t,\rho;s,\eta;c,H,b,m^2) \,  \mathcal{R}[|v|^p](s,\eta)\, \mathrm{d}\eta \,\mathrm{d}s.
\end{align}

The next step is to estimate the kernel function $E$ in the right-hand side of the last inequality on the shrunk $\eta$-interval of integration. First of all, from the Taylor expansion of the hypergeometric function 
\begin{align*}
\mathrm{F}(\tfrac 12 -\nu,\tfrac 12 -\nu; 1;\zeta)=\sum_{k=0}^\infty \frac{(\tfrac 12 -\nu)_k^2}{(1)_k \, k!}\zeta^k \qquad \mbox{for} \ \zeta\in [0,1),
\end{align*} where $(a)_0\doteq 1$ and $(a)_k\doteq a(a+1)\cdots (a+k-1)$ denotes the so-called Pochhammer symbol, we see immediately that we can estimate from below the factor involving the hypergeometric function in  $E(t,\rho;s,\eta;c,H,b,m^2)$ by the constant function $1$. Furthermore, since the two  exponential terms in \eqref{def E} are independent of $\eta$, the only factor that we actually have to estimate from below for $\eta\in [-(R+A(s)),R+A(s)]$ is $((\tfrac{c}{H}(\mathrm{e}^{Ht}+\mathrm{e}^{Hs}))^2-(\rho-\eta)^2)^{-\frac{1}{2}+\nu}$. 
Notice that we have to proceed in a different way in order to get such lower bound estimate depending on whether $\nu$ is smaller or greater than $ 1/2$.

Hereafter, for the sake of brevity, we use the notation $\phi(t)\doteq \frac{c}{H}\mathrm{e}^{Ht}$. In particular, we may express the amplitude of the forward light-cone as follows $A(t)=\phi(t)-\phi(0)$.

Let us begin with the case $\nu \in [0,\tfrac 12]$. Let us prove that in this case the following upper bound estimate holds
\begin{align} \label{upper bound phi(t)+phi(s)-rho+eta}
\phi(t)+\phi(s)-\rho+\eta\leqslant 2(\phi(t)-\rho)
\end{align}
for $s\in[0,s_0]$ and  $\eta\in [-(R+A(s)),R+A(s)]$.

 Clearly, \eqref{upper bound phi(t)+phi(s)-rho+eta} is equivalent to require that $\phi(t)-\rho-\eta\geqslant \phi(s)$. We check the validity of this inequality for $s\in[0,s_0]$ and $\eta\in [-(R+A(s)),R+A(s)]$ through a chain of inequalities
\begin{align*}
\phi(t)-\rho-\eta & \geqslant \phi(t)-\rho -R-A(s) \geqslant \phi(t)-\rho -R-A(s_0) \\
& = A(t)+ \phi(0)-\rho -R -\tfrac{1}{2}(A(t)-|\rho|-R) \geqslant \tfrac{1}{2}(A(t)-|\rho|-R) + \phi(0) \\ 
& \geqslant  A(s_0)+  \phi(0) \geqslant   A(s)+  \phi(0) = \phi(s),
\end{align*} where we used twice the condition $A(s)\leqslant A(s_0)$ and the identity $2A(s_0)=A(t)-|\rho|-R$.

In a completely analogous way, one proves that
\begin{align*} 
\phi(t)+\phi(s)-\rho+\eta\leqslant 2(\phi(t)+\rho)
\end{align*}
for $s\in[0,s_0]$ and  $\eta\in [-(R+A(s)),R+A(s)]$.

Therefore, combining \eqref{upper bound phi(t)+phi(s)-rho+eta} and the last inequality, when $\nu\leqslant \tfrac 12$ we can estimate
\begin{align}\label{important factor E case nu < 1/2}
((\phi(t)+\phi(s))^2-(\rho-\eta)^2)^{-\frac{1}{2}+\nu}\geqslant 2^{-1+2\nu} (\phi^2(t)-\rho^2)^{-\frac{1}{2}+\nu}
\end{align} for $s\in[0,s_0]$ and  $\eta\in [-(R+A(s)),R+A(s)]$.

On the other hand, for $\nu \geqslant \tfrac 12$, from the lower bound estimates
\begin{equation}\label{intermediate lb}
\begin{split}
\phi(t)+\phi(s)-\rho+\eta & \geqslant \phi(t)+\phi(s)-\rho-R -A(s) \geqslant \phi(t)-\rho+ \phi(0)-R ,\\
\phi(t)+\phi(s)+\rho-\eta & \geqslant \phi(t)+\phi(s)+\rho-R -A(s) \geqslant \phi(t)+\rho+ \phi(0)-R ,
\end{split}
\end{equation} for $s\in[0,s_0]$ and  $\eta\in [-(R+A(s)),R+A(s)]$, it follows that
\begin{align}\label{important factor E case nu > 1/2}
((\phi(t)+\phi(s))^2-(\rho-\eta)^2)^{-\frac{1}{2}+\nu}\geqslant ((\phi(t)+\phi(0)-R)^2- \rho^2)^{-\frac{1}{2}+\nu}. 
\end{align}

Combining \eqref{important factor E case nu < 1/2} and \eqref{important factor E case nu > 1/2}, we conclude that the kernel function in \eqref{second lower bound Rv} can be estimate from below in the following way
\begin{align}\label{kernel E lower bound}
 E(t,\rho;s,\eta;c,H,b,m^2) \gtrsim \mathrm{e}^{-(\frac{b}{2}+\nu H)t}  	\, \mathrm{e}^{(\frac{b}{2}-\nu H)s} ((\phi(t)+R_1)^2- \rho^2)^{-\frac{1}{2}+\nu}
\end{align} for $s\in[0,s_0]$ and  $\eta\in [-(R+A(s)),R+A(s)]$, where
\begin{align}\label{def R1}
R_1\doteq \begin{cases} 0 & \mbox{if}  \ \nu\leqslant \tfrac 12,\\ \phi(0)-R & \mbox{if} \ \nu > \tfrac 12. \end{cases}
\end{align}

\begin{remark} Notice that for $R\leqslant \phi(0)$, in \eqref{kernel E lower bound} we might consider $R_1=0$ even for $\nu>\tfrac 12$.
\end{remark}

\begin{remark} In the previous considerations we estimate from below the hypergeometric function by a constant. In the limit case $b^2=4m^2$ (that is, for $\nu=0$), we might think to employ the asymptotic estimate $\mathrm{F}(\frac{1}{2},\frac{1}{2};1;\zeta)\sim - \ln (1-\zeta)$ as $\zeta\to 1^-$ in order to improve this lower bound estimate. However, for $s\in[0,s_0]$ and $\eta\in [-(R+A(s)),R+A(s)]$, setting $$\zeta=\zeta(t,\rho;s,\eta;c,H)\doteq \frac{\left(\tfrac{c}{H}(\mathrm{e}^{Ht}-\mathrm{e}^{Hs})\right)^2-(\rho-\eta)^2}{\left(\tfrac{c}{H}(\mathrm{e}^{Ht}+\mathrm{e}^{Hs})\right)^2-(\rho-\eta)^2},$$  we have
\begin{align*}
-\ln(1-\zeta) & =\ln \left(\frac{(\phi(t)+\phi(s))^2-(\rho-\eta)^2}{4 \phi(t)\phi(s)}\right)\geqslant \ln \left(\frac{(\phi(t)+\phi(0)-R)^2-\rho^2}{4 \phi(t)\phi(s_0)}\right)\\ &= \ln \left(\frac{\phi(t)+\phi(0)-R+\rho}{2 \phi(t)}\right),
\end{align*} where  we used \eqref{intermediate lb} and $\phi(s_0)=\frac{1}{2}(\phi(t)+\phi(0)-\rho-R)$. Therefore, $-\ln(1-\zeta)$ does not provide an improvement in the lower bound estimate, since for large $t$ the argument of the logarithmic term on the right-hand side of the last inequality can be only estimated by a constant for $\rho \in [0,A(t)-R]$ (this is the actual range that we will consider for $\rho$ at the end of the present subsection).
\end{remark}
Now we plug the lower bound estimate from \eqref{kernel E lower bound} in \eqref{second lower bound Rv}. For $|\rho|\leqslant A(t)-R$ we have
\begin{align}\label{Radon transf 1st fund ineq}
\mathcal{R}[v](t,\rho) & \geqslant \mathrm{e}^{-(\frac{b}{2}+\nu H)t}  	\, ((\phi(t)+R_1)^2- \rho^2)^{-\frac{1}{2}+\nu}  \int_{0}^{s_0} \Gamma(s) \, \mathrm{e}^{(\frac{b}{2}-\nu H)s} \|v(s,\cdot)\|^{(\beta+1) p}_{L^p(\mathbb{R}^n)} \,\mathrm{d}s, 
\end{align} where we used the support condition for $\mathcal{R}[|v|^p]$ and Fubini's theorem to obtain
\begin{align*}
\int_{-(R+A(s))}^{R+A(s)}   \mathcal{R}[|v|^p](s,\eta)\, \mathrm{d}\eta \,\mathrm{d}s =\int_{\mathbb{R}}   \mathcal{R}[|v|^p](s,\eta)\, \mathrm{d}\eta \,\mathrm{d}s =  \|v(s,\cdot)\|^{p}_{L^p(\mathbb{R}^n)} .
\end{align*}

The inequality in \eqref{Radon transf 1st fund ineq} is the first crucial estimate to obtain the iteration frame for $ \|v(t,\cdot)\|^{p}_{L^p(\mathbb{R}^n)}$. Next step is to determine a lower bound for $ \|v(t,\cdot)\|^{p}_{L^p(\mathbb{R}^n)}$ with $\mathcal{R}[v]$ appearing in a nonlinear term on the right-hand side.

In order to derive this inequality, we will follow the approach from \cite[Section 5]{PalRei19}. We introduce the operator
\begin{align*}
\mathcal{T}(h)(\tau)\doteq |A(t)+R-\tau|^{-\frac{n-1}{2}}\int_\tau^{A(t)+R}h(r)|r-\tau|^{\frac{n-3}{2}}\mathrm{d}r
\end{align*} for any $\tau\in\mathbb{R}$ and any $h\in L^p(\mathbb{R})$.

In \cite[Section 5]{PalRei19} it is proved that $\mathcal{T}\in \mathcal{L}(L^p(\mathbb{R})\to L^p(\mathbb{R}))$ for any $p\in (1,\infty)$ and $n\geqslant 2$. Even though the function $A(t)$ in \cite{PalRei19} is a polynomial function (more precisely, $A(t)=\tfrac{1}{\ell+1}(t^{\ell+1}-1)$ for some $\ell\geqslant 0$), the proof of this result is actually independent of the explicit expression of $A(t)$ and it can be repeated verbatim in our case with $A(t)=\frac{c}{H}(\mathrm{e}^{H t}-1)$. 

We consider now the function
\begin{align*}
h(t,r)\doteq \begin{cases} |v(t,r)| r^{\frac{n-1}{p}} & \mbox{if} \ r\geqslant 0,\\ 0 & \mbox{if} \ r <0. \end{cases}
\end{align*} By  the boundedness of the operator $\mathcal{T}$ on $L^p(\mathbb{R})$, we have that $\|\mathcal{T}(h)(t,\cdot)\|_{L^p(\mathbb{R})}\lesssim \|h(t,\cdot)\|_{L^p(\mathbb{R})}$ holds uniformly with respect to $t\in [0,T)$.

Therefore,
\begin{align}
\int_{\mathbb{R}^n}|v(t,x)|^p\mathrm{d}x & =\omega_n \int_0^\infty |v(t,r)|^p r^{n-1} \mathrm{d}r=\omega_n  \| h(t,\cdot)\|^p_{L^p(\mathbb{R})} \notag \\ & \gtrsim \|\mathcal{T}(h)(t,\cdot)\|_{L^p(\mathbb{R})}^p =\int_{\mathbb{R}}|\mathcal{T}(h)|(t,\rho)|^p\mathrm{d}\rho \notag \\ &= \int_{\mathbb{R}}|A(t)+R-\rho|^{-\frac{n-1}{2}p}\left|\int_{\rho}^{A(t)+R} |v(t,r)|r^{\frac{n-1}{p}}|r-\rho|^{\frac{n-3}{2}}\mathrm{d}r\right|^p \mathrm{d}\rho \notag
\\
&\geqslant  \int_{0}^{A(t)+R}(A(t)+R-\rho)^{-\frac{n-1}{2}p}\left(\int_{\rho}^{A(t)+R} |v(t,r)|r^{\frac{n-1}{p}}(r-\rho)^{\frac{n-3}{2}}\mathrm{d}r\right)^p \mathrm{d}\rho. \label{lb int |v|^p oper mathcalT}
\end{align}

We have seen that $\mathcal{R}[v]$ is a nonnegative function by using an explicit integral representation. Moreover, by using the monotonicity of $\mathcal{R}$ and \eqref{Radon transform radial case}, we get
\begin{align}
0\leqslant \mathcal{R}[v](t,\rho) \leqslant \mathcal{R}[|v|](t,\rho) & =\omega_{n-1}\int_{|\rho|}^{A(t)+R}|v(t,r)|(r-|\rho|)^{\frac{n-3}{2}} (r+|\rho|)^{\frac{n-3}{2}} r \, \mathrm{d}r \notag  \\
& \leqslant 2^{\frac{n-3}{2}}\omega_{n-1}\int_{|\rho|}^{A(t)+R}|v(t,r)|(r-|\rho|)^{\frac{n-3}{2}} r^{\frac{n-1}{2}} \mathrm{d}r. \label{Radon transf upper bound}
\end{align}
Clearly, for $\rho\in [0,A(t)+R]$ and $r\in[\rho, A(t)+R]$ it results
\begin{align*}
r^{\frac{n-1}{p}}\geqslant (A(t)+R)^{-(n-1)[\frac 1p -\frac 12]_-} r^{\frac{n-1}{2}}\rho^{(n-1)[\frac 1p -\frac 12]_+}, 
\end{align*} where $[\tfrac 1p -\tfrac 12]_{\pm}$ denote the positive and negative part of $\tfrac 1p -\tfrac 12$, respectively. 

Combining \eqref{lb int |v|^p oper mathcalT}, \eqref{Radon transf upper bound} and the above inequality, we have
\begin{align}
\|v(t,\cdot)\|^p_{L^p(\mathbb{R}^n)} & \gtrsim (A(t)+R)^{-(n-1)[1 -\frac p2]_-}    \notag \\
& \qquad \times \int_{0}^{A(t)+R}\frac{\rho^{(n-1)[1 -\frac p2]_+} }{(A(t)+R-\rho)^{\frac{n-1}{2}p}}  \left(\int_{\rho}^{A(t)+R} |v(t,r)|r^{\frac{n-1}{2}}(r-\rho)^{\frac{n-3}{2}}\mathrm{d}r\right)^p \mathrm{d}\rho \notag \\
 & \gtrsim (A(t)+R)^{-(n-1)[1 -\frac p2]_-}   \int_{0}^{A(t)+R}\frac{\rho^{(n-1)[1 -\frac p2]_+} }{(A(t)+R-\rho)^{\frac{n-1}{2}p}}\big(\mathcal{R}[v](t,\rho)\big)^p \mathrm{d}\rho. \label{Radon transf 2nd fund ineq}
\end{align}

Finally, from \eqref{Radon transf 1st fund ineq} and \eqref{Radon transf 2nd fund ineq}, we obtain for $t\geqslant A^{-1}(R)$ the desired iteration frame 
\begin{align}
\|v(t,\cdot)\|^p_{L^p(\mathbb{R}^n)} & \geqslant   \frac{K \mathrm{e}^{-(\frac{b}{2}+\nu H)pt}}{ (A(t)+R)^{(n-1)[1 -\frac p2]_-}}   \int_{0}^{A(t)-R}\frac{\rho^{(n-1)[1 -\frac p2]_+} }{(A(t)+R-\rho)^{\frac{n-1}{2}p}} ((\phi(t)+R_1)^2- \rho^2)^{(-\frac{1}{2}+\nu)p} \notag \\
& \qquad \qquad \times \left(    \int_{0}^{A^{-1}(\frac 12 (A(t)-\rho-R))} \Gamma(s) \, \mathrm{e}^{(\frac{b}{2}-\nu H)s} \|v(s,\cdot)\|^{(\beta+1) p}_{L^p(\mathbb{R}^n)} \,\mathrm{d}s \right)^p \mathrm{d}\rho
\label{Iteration fram ||v(t)||^p Lp}
\end{align} for a suitable positive multiplicative constant $K=K(n,H,b,m^2)$. Needless to say, in \eqref{Iteration fram ||v(t)||^p Lp} only one between the factors $ (A(t)+R)^{-(n-1)[1 -\frac p2]_-}$ and $\rho^{(n-1)[1 -\frac p2]_+}$ is actually present for $p\neq 2$. Nevertheless, we will do the computations formally as if both were present in order to consider simultaneously the cases $p\in(1,2)$ and $p\geqslant 2$.

\begin{remark} \label{Remark n=1}
Let us underline explicitly that \eqref{Iteration fram ||v(t)||^p Lp} is true also for $n=1$. First, \eqref{Radon transf 1st fund ineq} can be obtained exactly as we did for $n\geqslant 2$, with the only difference that Yagdjian integral representation formula is applied now directly to $v$, that is,
\begin{align}\label{Radon transf 1st fund ineq n=1}
v(t,\rho) & \geqslant \mathrm{e}^{-(\frac{b}{2}+\nu H)t}  	\, ((\phi(t)+R_1)^2- \rho^2)^{-\frac{1}{2}+\nu}  \int_{0}^{s_0} \Gamma(s) \, \mathrm{e}^{(\frac{b}{2}-\nu H)s} \|v(s,\cdot)\|^{(\beta+1) p}_{L^p(\mathbb{R})} \,\mathrm{d}s.
\end{align}
On the other hand,  for $n=1$ \eqref{Radon transf 2nd fund ineq} can be replaced by the trivial inequality
\begin{align} \label{Radon transf 2nd fund ineq n=1}
\|v(t,\cdot)\|^p_{L^p(\mathbb{R})} & \geqslant \int_{0}^{A(t)+R}|v(t,\rho)|^p \mathrm{d}\rho.
\end{align} Hence, combining \eqref{Radon transf 1st fund ineq n=1} and \eqref{Radon transf 2nd fund ineq n=1}, we conclude the validity of \eqref{Iteration fram ||v(t)||^p Lp} for $n=1$ too.
\end{remark}

\subsection{Iteration argument  for $\|v(t,\cdot)\|^p_{L^p(\mathbb{R}^n)}$} \label{Subsection sequence lb Lp norm v(t)}

Our next goal is to derive a sequence of lower bound estimates for $\|v(t,\cdot)\|^p_{L^p(\mathbb{R}^n)}$  through the iteration argument \eqref{Iteration fram ||v(t)||^p Lp}.

 The starting point of our iteration procedure is given by \eqref{||v(t)||^p Lp lower bound step -1}.

Let us derive now a first lower bound estimate for $\|v(t,\cdot)\|^p_{L^p(\mathbb{R}^n)}$ with an additional polynomial growing factor.
Plugging \eqref{||v(t)||^p Lp lower bound step -1} in \eqref{Iteration fram ||v(t)||^p Lp}, for $t\geqslant A^{-1}(R)$ we get
\begin{align} 
\|v(t,\cdot)\|^p_{L^p(\mathbb{R}^n)} & \geqslant   K \widetilde{B}^q \varepsilon^{pq} \mathrm{e}^{-(\frac{b}{2}+\nu H)pt} (A(t)+R)^{-(n-1)[1 -\frac p2]_-}   \notag \\
& \qquad \times \int_{0}^{A(t)-R}\rho^{(n-1)[1 -\frac p2]_+} \frac{((\phi(t)+R_1)^2- \rho^2)^{(-\frac{1}{2}+\nu)p}}{(A(t)+R-\rho)^{\frac{n-1}{2}p}}  (I_0(t,\rho))^p \mathrm{d}\rho, \label{||v(t)||^p Lp lower bound step 0 intermediate}
\end{align} where, henceforth, $q\doteq (\beta+1)p$ and
\begin{align*}
I_0(t,\rho)\doteq  \int_{0}^{A^{-1}(\frac 12 (A(t)-\rho-R))} \Gamma(s) \, \mathrm{e}^{\left(\frac{b}{2}-\nu H-\frac{1}{2}(b+H)q+(n-1)H(\beta+1) -\frac{(n-1)}{2}Hq\right)s} \,\mathrm{d}s.
\end{align*}
We recall that $\Gamma(s)=\mu (1+s)^{\varsigma}\mathrm{e}^{\varrho_{\mathrm{crit}} s}$, being $\varrho_{\mathrm{crit}}$ defined by \eqref{def rho crit}. Therefore, using the actual value of $\Gamma(s)$, we have
\begin{align*}
I_0(t,\rho) = \int_{0}^{A^{-1}(\frac 12 (A(t)-\rho-R))} \mu \, (1+s)^\varsigma \, \mathrm{e}^{(\frac{n}{2}-\nu-\frac{1}{p})Hs} \,\mathrm{d}s,
\end{align*} where we used
\begin{align}
 & \varrho_{\mathrm{crit}}+\tfrac{b}{2}-\nu H-\tfrac{1}{2}(b+H)q+(n-1)H(\beta+1) -\tfrac{(n-1)}{2}Hq \notag \\
 & \qquad = \tfrac{1}{2}(b+nH)(q-1)+nH-\tfrac{H}{p}+ \tfrac{b}{2}-\nu H-\tfrac{1}{2}(b+H)q -\tfrac{(n-1)}{2}Hq \notag \\
 & \qquad = \tfrac{nH}{2}-\nu H-\tfrac{H}{p} \label{coefficient exp factor s integral}
\end{align} for the coefficient in the exponential term.

Hereafter, we consider only the case $\varsigma\leqslant 0$. The complementary case $\varsigma>0$ will be discuss at the end of Section \ref{Subsection anti dS nlin lb poly growth} (cf. Subsection \ref{Subsubsection varsigma >0}).

Therefore,
\begin{align*}
I_0(t,\rho) & \geqslant  \mu\,  (1+t)^\varsigma \int_{0}^{A^{-1}(\frac 12 (A(t)-\rho-R))}  \mathrm{e}^{(\frac{n}{2}-\nu-\frac{1}{p})Hs} \,\mathrm{d}s \\
& = \frac{ \mu}{(\frac{n}{2}-\nu-\frac{1}{p})H}  (1+t)^\varsigma   \left(\mathrm{e}^{(\frac{n}{2}-\nu-\frac{1}{p})HA^{-1}\left(\frac 12 (A(t)-\rho-R)\right)} -1\right) \\
& = \frac{ \mu}{(\frac{n}{2}-\nu-\frac{1}{p})H}  (1+t)^\varsigma   \left(\left(\tfrac{H}{2c}\right)^{\frac{n}{2}-\nu-\frac{1}{p}}\left(A(t)-\rho-R+\tfrac{2c}{H}\right)^{\frac{n}{2}-\nu-\frac{1}{p}} -1\right),
\end{align*} 
where we used the analytic expression of the inverse function of $A$ given by 
\begin{align}\label{A inverse function}
A^{-1}(\sigma)= \frac{1}{H}\ln \left(\frac{H\sigma}{c}+1\right).
\end{align}

Let $b_0>0$ be a fixed parameter. For $\rho\leqslant A(t)-R-b_0\frac{c}{H}$ we may estimate $I_0(t,\rho)$ from below as follows
\begin{align*}
I_0(t,\rho) & \geqslant\widehat{B} (1+t)^\varsigma  \left(A(t)-\rho-R+\tfrac{2c}{H}\right)^{\frac{n}{2}-\nu-\frac{1}{p}}, 
\end{align*} where $\widehat{B}\doteq \mu (\frac{n}{2}-\nu-\frac{1}{p})^{-1}H^{-1} \left(1-(\tfrac{b_0}{2}+1)^{-\frac{n}{2}+\nu+\frac{1}{p}}\right)\left(\frac{H}{2c}\right)^{\frac{n}{2}-\nu-\frac{1}{p}}$. Plugging the last lower bound for $I_0(t,\rho)$ in \eqref{||v(t)||^p Lp lower bound step 0 intermediate}, we get
\begin{align} \label{||v(t)||^p Lp lower bound step 0 intermediate 2}
\|v(t,\cdot)\|^p_{L^p(\mathbb{R}^n)} & \geqslant   K \widetilde{B}^q \widehat{B}^p \varepsilon^{pq} \mathrm{e}^{-(\frac{b}{2}+\nu H)pt} (A(t)+R)^{-(n-1)[1 -\frac p2]_-}    (1+t)^{\varsigma p} J_0(t)
\end{align}  for $t\geqslant A^{-1}(R+b_0\frac{c}{H})$, where
\begin{align*}
J_0(t)\doteq \int_{0}^{A(t)-R-b_0\frac{c}{H}}\rho^{(n-1)[1 -\frac p2]_+}    \frac{((\phi(t)+R_1)^2- \rho^2)^{(-\frac{1}{2}+\nu)p}}{(A(t)+R-\rho)^{\frac{n-1}{2}p}} \left(A(t)-\rho-R+\tfrac{2c}{H}\right)^{\frac{n}{2} p-\nu p-1}  \mathrm{d}\rho.
\end{align*}
The next step consists in estimating from below the integral $J_0(t)$.

 First, we consider the factor $((\phi(t)+R_1)^2- \rho^2)^{(-\frac{1}{2}+\nu)p}$. Recalling that the value of $R_1$ depends on the range for $\nu$ (cf. \eqref{def R1} for the definition of $R_1$), we derive a lower bound for this factor separately in the case $\nu\leqslant \frac{1}{2}$ and in the case $\nu > \frac{1}{2}$. 

For $\nu\leqslant \frac{1}{2}$, since the power $(-\frac{1}{2}+\nu)p$ is nonpositive, we consider an upper bound for $(\phi(t)+R_1)^2- \rho^2$. For $\rho\in [0,A(t)-R]$, since $R_1=0$ in this case, we have
\begin{align*}
\phi(t)+R_1+ \rho &\leqslant \phi(t)+A(t)-R=2\phi (t) -\tfrac cH -R\leqslant  2\phi (t),
\end{align*} and
\begin{align*}
\phi(t)+R_1- \rho & = A(t)-\rho +\tfrac cH \leqslant \begin{cases} A(t)-\rho-R +\tfrac{2c}{H} & \mbox{if} \ R\leqslant \tfrac cH, \\ A (t)-\rho+R & \mbox{if} \ R >\tfrac cH. \end{cases} 
\end{align*} Thus, for $\nu\leqslant \frac{1}{2}$ and $\rho\in [0,A(t)-R]$ we obtained
\begin{align}
((\phi(t)+R_1)^2- \rho^2)^{(-\frac{1}{2}+\nu)p} \geqslant  \begin{cases}\left(\frac{2c}{H}\right)^{(-\frac{1}{2}+\nu)p} \mathrm{e}^{(-\frac{1}{2}+\nu)Hpt} (A(t)-\rho-R +\tfrac{2c}{H})^{(-\frac{1}{2}+\nu)p}  & \mbox{if} \ R\leqslant \tfrac cH, \\ \left(\frac{2c}{H}\right)^{(-\frac{1}{2}+\nu)p} \mathrm{e}^{(-\frac{1}{2}+\nu)Hpt} (A(t)-\rho+R)^{(-\frac{1}{2}+\nu)p}  & \mbox{if} \ R>\tfrac cH. \end{cases} \label{lower bound ((phi(t)+R1)^2-rho^2)^((-1/2+nu)p) case nu <1/2}
\end{align}
We consider now $\nu > \frac{1}{2}$. In this case we determine a lower bound for $(\phi(t)+R_1)^2- \rho^2$.  For $\rho\in [0,A(t)-R]$, since $R_1=\frac{c}{H}-R$ in this case, we have
\begin{align*}
\phi(t)+R_1+ \rho & \geqslant \phi(t)+\tfrac{c}{H}-R \geqslant \begin{cases} \phi(t) & \mbox{if} \ R\leqslant \tfrac cH, \\ \frac{1}{2} \phi(t) & \mbox{if} \ R > \tfrac cH \ \mbox{and} \ t\geqslant A^{-1}(2R-\frac{3c}{H}), \end{cases}  
\end{align*} and $\phi(t)+R_1- \rho  = A(t)-\rho-R+\tfrac{2c}{H}$.

 We emphasize that the lower bound for $\|v(t,\cdot)\|_{L^p(\mathbb{R}^n)}^p$ that we are going to prove will be valid for $t\geqslant A^{-1}(a_0 R+b_0 \frac{c}{H})$ and for suitable $a_0\geqslant 2$ and $b_0>0$. In particular, we will use the inequality $\phi(t)+R_1+ \rho \geqslant \frac{1}{2} \phi(t)$ when  $ R > \tfrac cH $ without specifying the further condition on $t$.

Hence,  for $\nu>\frac{1}{2}$ and $\rho\in [0,A(t)-R]$ we proved that
\begin{align}
((\phi(t)+R_1)^2- \rho^2)^{(-\frac{1}{2}+\nu)p} \geqslant  \begin{cases}\left(\frac{c}{H}\right)^{(-\frac{1}{2}+\nu)p} \mathrm{e}^{(-\frac{1}{2}+\nu)Hpt} (A(t)-\rho-R +\tfrac{2c}{H})^{(-\frac{1}{2}+\nu)p}  & \mbox{if} \ R\leqslant \tfrac cH, \\ \left(\frac{c}{2H}\right)^{(-\frac{1}{2}+\nu)p} \mathrm{e}^{(-\frac{1}{2}+\nu)Hpt} (A(t)-\rho-R +\tfrac{2c}{H})^{(-\frac{1}{2}+\nu)p}  & \mbox{if} \ R> \tfrac cH. \end{cases}  \label{lower bound ((phi(t)+R1)^2-rho^2)^((-1/2+nu)p) case nu >1/2}
\end{align}

Next, we consider the factor $(A(t)+R-\rho)^{-\frac{n-1}{2}p}$ in $J_0(t)$. Notice that in the case $\nu\leqslant \frac{1}{2}$ when $R> \frac{c}{H}$, from \eqref{lower bound ((phi(t)+R1)^2-rho^2)^((-1/2+nu)p) case nu <1/2} we see that we have to consider actually the factor $(A(t)+R-\rho)^{(-\frac{n}{2}+\nu)p}$. From \eqref{1st lb for V from V0 is dominant} it follows that $(-\frac{n}{2}+\nu)p<0$. Hence, in both cases we are interested in an upper bound estimate for the term $A(t)+R-\rho$. By straightforward computations, we get that 
\begin{align}\label{A(t)-rho+R upper bound}
A(t)-\rho-R+\tfrac{2c}{H}\geqslant \tfrac{a_0-1}{a_0+1} (A(t)-\rho +R) \qquad \mbox{for any} \ \ \rho\leqslant A(t)-a_0 R+ (a_0+1)\tfrac{c}{H}. 
\end{align} Combining \eqref{lower bound ((phi(t)+R1)^2-rho^2)^((-1/2+nu)p) case nu <1/2}, \eqref{lower bound ((phi(t)+R1)^2-rho^2)^((-1/2+nu)p) case nu >1/2} and \eqref{A(t)-rho+R upper bound} and shrinking the domain of integration in $J_0(t)$, we get
\begin{align*}
J_0(t)\geqslant \bar{B} \mathrm{e}^{(-\frac{1}{2}+\nu)Hpt} \int_0^{A(t)-a_0R-b_0\frac{c}{H}} \rho^{(n-1)[1 -\frac p2]_+}   \left(A(t)-\rho-R+\tfrac{2c}{H}\right)^{-1} \mathrm{d}\rho
\end{align*} for $t\geqslant A^{-1}(a_0R+b_0 \frac{c}{H})$, where 
\begin{align}
\bar{B} \doteq \begin{cases}
\left(	\frac{2c}{H}\right)^{(-\frac{1}{2}+\nu)p} \left(\frac{a_0-1}{a_0+1}\right)^{\frac{n-1}{2}p} & \mbox{if} \ \nu\leqslant \frac{1}{2} \ \mbox{and} \ R\leqslant \frac{c}{H}, \\
\left(	\frac{2c}{H}\right)^{(-\frac{1}{2}+\nu)p} \left(\frac{a_0-1}{a_0+1}\right)^{(\frac{n}{2}-\nu)p} & \mbox{if} \ \nu\leqslant \frac{1}{2} \ \mbox{and} \ R>\frac{c}{H},\\
\left(	\frac{c}{H}\right)^{(-\frac{1}{2}+\nu)p} \left(\frac{a_0-1}{a_0+1}\right)^{\frac{n-1}{2}p} & \mbox{if} \ \nu > \frac{1}{2} \ \mbox{and} \ R\leqslant \frac{c}{H},\\
\left(	\frac{c}{2H}\right)^{(-\frac{1}{2}+\nu)p} \left(\frac{a_0-1}{a_0+1}\right)^{\frac{n-1}{2}p} & \mbox{if} \ \nu > \frac{1}{2} \ \mbox{and} \ R> \frac{c}{H}.
\end{cases} \label{B bar def}
\end{align}
Then, we shrink further the domain of integration in the right-hand side of the last inequality by increasing the bottom of the interval of integration from $0$ to $\delta({A(t)-a_0R-b_0\frac{c}{H}})$ for some $\delta\in (0,1)$, obtaining for $t\geqslant A^{-1}(a_0R+b_0 \frac{c}{H})$
\begin{align}
J_0(t) &\geqslant \bar{B} \mathrm{e}^{(-\frac{1}{2}+\nu)Hpt} \int_{\delta({A(t)-a_0R-b_0\frac{c}{H}})}^{A(t)-a_0R-b_0\frac{c}{H}} \rho^{(n-1)[1 -\frac p2]_+}   \left(A(t)-\rho-R+\tfrac{2c}{H}\right)^{-1} \mathrm{d}\rho \notag \\
&\geqslant \frac{\bar{B} \mathrm{e}^{(-\frac{1}{2}+\nu)Hpt}}{\delta^{-(n-1)[1 -\frac p2]_+}}  (A(t)-a_0R-b_0\tfrac{c}{H})^{(n-1)[1 -\frac p2]_+} \int_{\delta({A(t)-a_0R-b_0\frac{c}{H}})}^{A(t)-a_0R-b_0\frac{c}{H}}   \frac{ \mathrm{d}\rho }{A(t)-\rho-R+\tfrac{2c}{H}}\notag \\
& = \frac{\bar{B} \mathrm{e}^{(-\frac{1}{2}+\nu)Hpt}}{\delta^{-(n-1)[1 -\frac p2]_+}}  (A(t)-a_0R-b_0\tfrac{c}{H})^{(n-1)[1 -\frac p2]_+} \ln \left(\frac{A(t)+\frac{\delta a_0-1}{1-\delta}R+\frac{\delta b_0+2}{1-\delta}\frac{c}{H}}{\frac{a_0-1}{1-\delta}R+\frac{b_0+2}{1-\delta}\frac{c}{H}}\right). \label{final lb for J0}
\end{align} We emphasize that a more restrictive range for $\delta$ is going to be prescribed in the inductive step.

Finally, plugging \eqref{final lb for J0} in \eqref{||v(t)||^p Lp lower bound step 0 intermediate 2}, for  $t\geqslant A^{-1}(a_0R+b_0 \frac{c}{H})$ we get
\begin{align}
\|v(t,\cdot)\|^p_{L^p(\mathbb{R}^n)} & \geqslant   B_0 \varepsilon^{pq} (1+t)^{\varsigma p}  \mathrm{e}^{-\frac{1}{2}(b+H)pt}  (A(t)+R)^{-(n-1)[1 -\frac p2]_-} \notag  \\ & \qquad \times (A(t)-a_0R-b_0\tfrac{c}{H})^{(n-1)[1 -\frac p2]_+} \ln \left(\frac{A(t)+\frac{\delta a_0-1}{1-\delta}R+\frac{\delta b_0+2}{1-\delta}\frac{c}{H}}{\frac{a_0-1}{1-\delta}R+\frac{b_0+2}{1-\delta}\frac{c}{H}}\right), \label{||v(t)||^p Lp lower bound step 0}
\end{align} where $B_0=B_0(n,c,H,b,m^2,p,\mu,v_0,v_1,R,a_0,b_0,\delta) \doteq \delta^{(n-1)[1 -\frac p2]_+} K \widetilde{B}^q \widehat{B}^p \bar{B}$. We recall that the only assumptions on the parameters $a_0$ and $b_0$ that we did in order to obtain \eqref{||v(t)||^p Lp lower bound step 0} are the following
\begin{align}\label{condition a0, b0 first step}
a_0\geqslant 2 \quad \mbox{and} \quad b_0 >0.
\end{align} In the next subsection, however, we will require further conditions on $a_0$, cf. \eqref{condition a0 inductive step} and \eqref{condition a0 final lb V(t)}.

We stress that, since the amplitude function $A$ for the light-cone grows exponentially, the logarithmic term in \eqref{||v(t)||^p Lp lower bound step 0} provides, together with $(1+t)^{\varsigma p}$, a polynomially increasing factor. This factor constitutes the improvement with respect to the estimate in \eqref{||v(t)||^p Lp lower bound step -1} provided that $\varsigma \in (-\frac{1}{p},0]$.

The next step is to prove that $\|v(t,\cdot)\|^p_{L^p(\mathbb{R}^n)}$ satisfies the following sequence of lower bound estimates
\begin{align}
\|v(t,\cdot)\|^p_{L^p(\mathbb{R}^n)} & \geqslant   B_j \varepsilon^{pq^{j+1}} (1+t)^{\varsigma p \frac{q^{j+1}-1}{q-1}}  \mathrm{e}^{-\frac{1}{2}(b+H)pt}  (A(t)+R)^{-(n-1)[1 -\frac p2]_-} \notag  \\ & \quad \times (A(t)-a_jR-b_j\tfrac{c}{H})^{(n-1)[1 -\frac p2]_+} \left(\ln \! \left(\frac{A(t)+\frac{\delta a_j-1}{1-\delta}R+\frac{\delta b_j+2}{1-\delta}\frac{c}{H}}{\frac{a_j-1}{1-\delta}R+\frac{b_j+2}{1-\delta}\frac{c}{H}}\right)\!\right)^{\frac{q^{j+1}-1}{q-1}}
\label{||v(t)||^p Lp lower bound step j}
\end{align}  for $t\geqslant A^{-1}(a_j R+b_j \frac{c}{H})$, where $\{B_j\}_{j\in\mathbb{N}}$ is a suitable sequence of positive real numbers that we will determine iteratively during the proof and
\begin{align}
a_j & \doteq (a_0-1) \left(\tfrac{4}{1-\delta}\right)^j+1, \label{def aj geometric sequence}\\
b_j & \doteq (b_0+2) \left(\tfrac{4}{1-\delta}\right)^j-2. \label{def bj geometric sequence}
\end{align} 
Clearly, we have already proved \eqref{||v(t)||^p Lp lower bound step j} for $j=0$, namely, \eqref{||v(t)||^p Lp lower bound step 0}.
We are going to prove \eqref{||v(t)||^p Lp lower bound step j} by induction. Let us assume that \eqref{||v(t)||^p Lp lower bound step j} holds for some $j$, with $j\geqslant 0$. We will prove \eqref{||v(t)||^p Lp lower bound step j} for $j+1$, determining the value of $B_{j+1}$ in terms of $B_{j}$. According to this goal, we plug \eqref{||v(t)||^p Lp lower bound step j} into the iteration frame in \eqref{Iteration fram ||v(t)||^p Lp}, obtaining 
\begin{align} 
\|v(t,\cdot)\|^p_{L^p(\mathbb{R}^n)} & \geqslant   K   B_j^q \varepsilon^{pq^{j+2}} \mathrm{e}^{-(\frac{b}{2}+\nu H)pt} (A(t)+R)^{-(n-1)[1 -\frac p2]_-}  J_{j+1}(t) ,
 \label{||v(t)||^p Lp lower bound step j intermediate}
\end{align} where
\begin{align}\label{def J j+i}
 J_{j+1}(t) & \doteq \int_{0}^{A(t)-(2a_j+1)R-2b_j \frac{c}{H}} \rho^{(n-1)[1 -\frac p2]_+} \frac{((\phi(t)+R_1)^2- \rho^2)^{(-\frac{1}{2}+\nu)p}}{(A(t)+R-\rho)^{\frac{n-1}{2}p}} (I_{j+1}(t,\rho))^p \mathrm{d}\rho
 \end{align} for $t\geqslant A^{-1}((2a_j+1)R+2b_j \tfrac{c}{H})$ and
 \begin{align} 
I_{j+1}(t,\rho) \doteq   & \int_{A^{-1}(a_jR+b_j \frac{c}{H})}^{A^{-1}(\frac 12 (A(t)-\rho-R))} \Gamma(s) \, (1+s)^{\varsigma  \frac{q^{j+2}-q}{q-1}}  \, \mathrm{e}^{(\frac{b}{2}-\nu H-\frac{b+H}{2}q)s}  \notag  \\ 
& \qquad \qquad \qquad \times \frac{(A(s)-a_jR-b_j\tfrac{c}{H})^{(n-1)(\beta+1)[1 -\frac p2]_+}}{(A(s)+R)^{(n-1)(\beta+1))[1 -\frac p2]_-} } \notag  \\
 & \qquad \qquad \qquad \times \left(\ln \! \left(\frac{A(s)+\frac{\delta a_j-1}{1-\delta}R+\frac{\delta b_j+2}{1-\delta}\frac{c}{H}}{\frac{a_j-1}{1-\delta}R+\frac{b_j+2}{1-\delta}\frac{c}{H}}\right)\!\right)^{(\beta+1)\frac{q^{j+1}-1}{q-1}} \mathrm{d}s. \label{def I j+i}
\end{align} Notice that in \eqref{def J j+i} we shrank the $\rho$-domain of integration in order to have a nonempty $s$-domain of integration when using \eqref{||v(t)||^p Lp lower bound step j} to obtain $I_{j+1}(t,\rho)$.

We begin by deriving a lower bound estimate for $I_{j+1}(t,\rho)$. By using \eqref{coefficient exp factor s integral}, we can rewrite the first three factors in $I_{j+1}(t,\rho)$ as follows:
\begin{align*}
 \Gamma(s) \, (1+s)^{\varsigma  \frac{q^{j+2}-q}{q-1}}  \, \mathrm{e}^{(\frac{b}{2}-\nu H-\frac{b+H}{2}q)s} & = \mu (1+s)^{\varsigma +\varsigma \frac{q^{j+2}-q}{q-1}}  \, \mathrm{e}^{(\varrho_{\mathrm{crit}}+\frac{b}{2}-\nu H-\frac{b+H}{2}q)s} \\
  & = \mu (1+s)^{\varsigma \frac{q^{j+2}-1}{q-1}}  \, \mathrm{e}^{(\frac{n}{2}-\nu-\frac{1}{p})Hs} \, \mathrm{e}^{-(n-1)H(\beta+1)(1-\frac{p}{2})s}.
\end{align*}
By straightforward computations, we find that $A(s)-a_jR-b_j\tfrac{c}{H}\geqslant \tfrac 12 \phi(s)$ if and only if $s\geqslant A^{-1}(2a_jR+(2b_j+1)\tfrac{c}{H})$, while $A(s)+R\leqslant 2 \phi(s)$ if and only if $s\geqslant A^{-1}(R-\tfrac{2c}{H})$. Therefore, shrinking the domain of integration in \eqref{def I j+i} to $[A^{-1}(2a_jR+(2b_j+1)\tfrac{c}{H}),A^{-1}(\tfrac 12 (A(t)-\rho-R))]$ and working with $\rho\in [0,A(t)-(4a_j+1)R-2(2b_j+1)\tfrac{c}{H}]$ in \eqref{def J j+i}, we may estimate
\begin{align*}
\frac{(A(s)-a_jR-b_j\tfrac{c}{H})^{(n-1)(\beta+1)[1 -\frac p2]_+}}{(A(s)+R)^{(n-1)(\beta+1))[1 -\frac p2]_-} } \geqslant 2^{-(n-1)(\beta+1)|1-\frac{p}{2}|}\Big(\tfrac{c}{H}\mathrm{e}^{Hs}\Big)^{(n-1)(\beta+1)(1-\frac{p}{2})}.
\end{align*}
Hence, for $\rho\in [0,A(t)-(4a_j+1)R-2(2b_j+1)\tfrac{c}{H}]$ we have
\begin{align} \label{lower bound I j+1 intermediate}
I_{j+1}(t,\rho) & \geqslant  2^{-(n-1)(\beta+1)|1-\frac{p}{2}|}\mu \left(\tfrac{c}{H}\right)^{(n-1)(\beta+1)(1-\frac{p}{2})} (1+t)^{\varsigma \frac{q^{j+2}-1}{q-1}}   \widetilde{I}_{j+1}(t,\rho)  , 
\end{align} where
\begin{align}\label{def tilde I j+1}
  \widetilde{I}_{j+1}(t,\rho) &   \doteq  \int_{A^{-1}(2a_jR+(2b_j+1) \frac{c}{H})}^{A^{-1}(\frac 12 (A(t)-\rho-R))} \! \mathrm{e}^{(\frac{n}{2}-\nu-\frac{1}{p})Hs} \! \left(\ln \! \left(\frac{A(s)+\frac{\delta a_j-1}{1-\delta}R+\frac{\delta b_j+2}{1-\delta}\frac{c}{H}}{\frac{a_j-1}{1-\delta}R+\frac{b_j+2}{1-\delta}\frac{c}{H}}\right)\!\right)^{(\beta+1)\frac{q^{j+1}-1}{q-1}} \mathrm{d}s.
\end{align}
In order to get a lower bound estimate for $\widetilde{I}_{j+1}(t,\rho)$ we shrink further the domain of integration with respect to $s$. 
 Since the inequality $A^{-1}(\frac 14 (A(t)-\rho-R))\geqslant A^{-1}(2a_jR+(2b_j+1) \frac{c}{H})$ holds for $\rho\in [0,A(t)-(8a_j+1)R-4(2b_j+1)\tfrac{c}{H}]$, we may reduce the domain of integration in \eqref{def tilde I j+1} to $\left[A^{-1}(\tfrac 14 (A(t)-\rho-R)),A^{-1}\left(\tfrac 12 (A(t)-\rho-R)\right)\right]$ for $\rho$ in this interval. 

Thus, for $\rho\in [0,A(t)-(8a_j+1)R-4(2b_j+1)\tfrac{c}{H}]$ we get
\begin{align}
 \widetilde{I}_{j+1}(t,\rho) &   \geqslant  \int_{A^{-1}(\frac 14 (A(t)-\rho-R))}^{A^{-1}(\frac 12 (A(t)-\rho-R))} \! \mathrm{e}^{(\frac{n}{2}-\nu-\frac{1}{p})Hs} \! \left(\ln \! \left(\frac{A(s)+\frac{\delta a_j-1}{1-\delta}R+\frac{\delta b_j+2}{1-\delta}\frac{c}{H}}{\frac{a_j-1}{1-\delta}R+\frac{b_j+2}{1-\delta}\frac{c}{H}}\right)\!\right)^{(\beta+1)\frac{q^{j+1}-1}{q-1}} \mathrm{d}s \notag \\
 &   \geqslant  \left(\ln \! \left(\frac{A(t)-\rho+\left(\frac{4(\delta a_j-1)}{1-\delta}-1\right)R+\frac{4(\delta b_j+2)}{1-\delta}\frac{c}{H}}{\frac{4(a_j-1)}{1-\delta}R+\frac{4(b_j+2)}{1-\delta}\frac{c}{H}}\right)\!\right)^{(\beta+1)\frac{q^{j+1}-1}{q-1}} \widehat{I}_{j+1}(t,\rho),
 \label{lower bound tilde I j+1}
\end{align} where
\begin{align*}
 \widehat{I}_{j+1}(t,\rho) &\doteq  \int_{A^{-1}(\frac 14 (A(t)-\rho-R))}^{A^{-1}(\frac 12 (A(t)-\rho-R))} \! \mathrm{e}^{(\frac{n}{2}-\nu-\frac{1}{p})Hs} \mathrm{d}s. 
\end{align*} 
We have
\begin{align}
\widehat{I}_{j+1}(t,\rho) 
 &   =   \frac{\left(A(t)-\rho-R+\tfrac{2c}{H}\right)^{\frac{n}{2}-\nu-\frac{1}{p}}}{(\frac{n}{2}-\nu-\frac{1}{p})H \left(\tfrac{2c}{H}\right)^{\frac{n}{2}-\nu-\frac{1}{p}}}  \left(1-\mathrm{e}^{(\frac{n}{2}-\nu-\frac{1}{p})H\left(A^{-1}\left(\frac 14 (A(t)-\rho-R)\right)-A^{-1}\left(\frac{1}{2} (A(t)-\rho-R)\right)\right)}\right), \label{lower bound hat I j+1}
\end{align} where we used \eqref{A inverse function} to get 
\begin{align*} 
\mathrm{e}^{(\frac{n}{2}-\nu-\frac{1}{p})HA^{-1}(\frac 12 (A(t)-\rho-R))} &= \left(\tfrac{H}{2c}\right)^{\frac{n}{2}-\nu-\frac{1}{p}}\left(A(t)-\rho-R+\tfrac{2c}{H}\right)^{\frac{n}{2}-\nu-\frac{1}{p}}. 
\end{align*} We estimate now the last factor in \eqref{lower bound hat I j+1} for  $\rho\in [0,A(t)-(8a_j+1)R-4(2b_j+1)\tfrac{c}{H}]$ as follows:
\begin{align} 
 & 1- \exp\left(\left(\tfrac{n}{2}-\nu-\tfrac{1}{p}\right)H\left(A^{-1}\left(\tfrac 14 (A(t)-\rho-R)\right)-A^{-1}\left(\tfrac{1}{2} (A(t)-\rho-R)\right)\right)\right)  \notag \\ & \quad = 1- \exp\left(\left(\tfrac{n}{2}-\nu-\tfrac{1}{p}\right)\left(\ln\left(\frac{\tfrac{H }{4c}(A(t)-\rho-R)+1}{\tfrac{H}{2c} (A(t)-\rho-R)+1}\right)\right)\right) \notag \\
   & \quad = 1- \exp\left(\left(\tfrac{n}{2}-\nu-\tfrac{1}{p}\right)\left(\ln\left(\frac 12 +\frac{\frac{c}{H}}{ A(t)-\rho-R+\tfrac{2c}{H}}\right)\right)\right) \notag \\
   & \quad \geqslant 1- \exp\left(\left(\tfrac{n}{2}-\nu-\tfrac{1}{p}\right)\left(\ln\left(\frac 12 +\frac{\frac{c}{H}}{ 8a_j R+ 2(4b_j+3)\tfrac{c}{H}}\right)\right)\right). \label{difference exponential factor}
\end{align} Let us denote by $d_j =d_j(n,c,H,\nu,p,R,a_0,b_0)$ the right-hand side of \eqref{difference exponential factor}. Since $\{d_j\}_{j\in \mathbb{N}}$ is an increasing and bounded sequence, we can simply estimate the left-hand side of \eqref{difference exponential factor} from below with $d_0$. 
Summarizing, we proved
\begin{align*}
\widehat{I}_{j+1}(t,\rho)  \geqslant \frac{d_0}{(\frac{n}{2}-\nu-\frac{1}{p})H \left(\tfrac{2c}{H}\right)^{\frac{n}{2}-\nu-\frac{1}{p}}} \left(A(t)-\rho-R+\tfrac{2c}{H}\right)^{\frac{n}{2}-\nu-\frac{1}{p}}.
\end{align*}
Hence, combining \eqref{lower bound tilde I j+1} and the last inequality, we conclude that 
\begin{align*}
\widetilde{I}_{j+1}(t,\rho) &  \geqslant \widetilde{C} \left(A(t)-\rho-R+\tfrac{2c}{H}\right)^{\frac{n}{2}-\nu-\frac{1}{p}}  \\ & \qquad \qquad \times  \left(\ln \! \left(\frac{A(t)-\rho+\left(\frac{4(\delta a_j-1)}{1-\delta}-1\right)R+\frac{4(\delta b_j+2)}{1-\delta}\frac{c}{H}}{\frac{4(a_j-1)}{1-\delta}R+\frac{4(b_j+2)}{1-\delta}\frac{c}{H}}\right)\!\right)^{(\beta+1)\frac{q^{j+1}-1}{q-1}} 
\end{align*} for $\rho\in [0,A(t)-(8a_j+1)R-4(2b_j+1)\tfrac{c}{H}]$, where the multiplicative constant is given by $\widetilde{C}\doteq \big(\frac{n}{2}-\nu-\frac{1}{p}\big)^{-1}H^{-1} \left(\tfrac{H}{2c}\right)^{\frac{n}{2}-\nu-\frac{1}{p}} d_0$.

Then, plugging the last lower bound for $\widetilde{I}_{j+1}(t,\rho)$ in \eqref{lower bound I j+1 intermediate}, we arrive at
\begin{align*} 
I_{j+1}(t,\rho) & \geqslant \widehat{C}  (1+t)^{\varsigma \frac{q^{j+2}-1}{q-1}}   \left(A(t)-\rho-R+\tfrac{2c}{H}\right)^{\frac{n}{2}-\nu-\frac{1}{p}}\\ & \qquad \qquad  \times \left(\ln \! \left(\frac{A(t)-\rho+\left(\frac{4(\delta a_j-1)}{1-\delta}-1\right)R+\frac{4(\delta b_j+2)}{1-\delta}\frac{c}{H}}{\frac{4(a_j-1)}{1-\delta}R+\frac{4(b_j+2)}{1-\delta}\frac{c}{H}}\right)\!\right)^{(\beta+1)\frac{q^{j+1}-1}{q-1}} 
\end{align*} for $\rho\in [0,A(t)-(8a_j+1)R-4(2b_j+1)\tfrac{c}{H}]$, where $\widehat{C}\doteq 2^{-(n-1)(\beta+1)|1-\frac{p}{2}|}\mu \left(\tfrac{c}{H}\right)^{(n-1)(\beta+1)(1-\frac{p}{2})} \widetilde{C}$.

Next, after shrinking the domain of integration with respect to $\rho$ as we described in the previous steps, we plug the obtained lower bound for $I_{j+1}(t,\rho)$ in \eqref{def J j+i}, obtaining
\begin{align*}
J_{j+1}(t) \geqslant \widehat{C}^p  (1+t)^{\varsigma p \frac{q^{j+2}-1}{q-1}} &  \int_{0}^{A(t)-(8a_j+1)R-4(2b_j+1) \frac{c}{H}} \rho^{(n-1)[1 -\frac p2]_+}  
\\ & \quad   \times \frac{((\phi(t)+R_1)^2- \rho^2)^{(-\frac{1}{2}+\nu)p}}{(A(t)+R-\rho)^{\frac{n-1}{2}p}} \left(A(t)-\rho-R+\tfrac{2c}{H}\right)^{(\frac{n}{2}-\nu)p-1} 
\\ & \quad  \times \left(\ln \! \left(\frac{A(t)-\rho+\left(\frac{4(\delta a_j-1)}{1-\delta}-1\right)R+\frac{4(\delta b_j+2)}{1-\delta}\frac{c}{H}}{\frac{4(a_j-1)}{1-\delta}R+\frac{4(b_j+2)}{1-\delta}\frac{c}{H}}\right)\!\right)^{\frac{q^{j+2}-q}{q-1}} \mathrm{d}\rho.
\end{align*} Since the sequences $\{a_j\}_{j\in\mathbb{N}}$, $\{b_j\}_{j\in\mathbb{N}}$ are strictly increasing, we may apply the same estimates for the factors in the middle line of the previous inequality that we used in the base case $j=0$, namely, \eqref{lower bound ((phi(t)+R1)^2-rho^2)^((-1/2+nu)p) case nu <1/2}, \eqref{lower bound ((phi(t)+R1)^2-rho^2)^((-1/2+nu)p) case nu >1/2} and \eqref{A(t)-rho+R upper bound}, arriving at 
\begin{align}
J_{j+1}(t) & \geqslant \bar{B} \widehat{C}^p  (1+t)^{\varsigma p \frac{q^{j+2}-1}{q-1}}  \mathrm{e}^{(-\frac{1}{2}+\nu)Hpt}   \int_{0}^{A(t)-(8a_j+1)R-4(2b_j+1) \frac{c}{H}} \frac{\rho^{(n-1)[1 -\frac p2]_+}}{  \left(A(t)-\rho-R+\tfrac{2c}{H}\right)} \notag
\\ & \qquad \qquad \times \left(\ln \! \left(\frac{A(t)-\rho+\left(\frac{4(\delta a_j-1)}{1-\delta}-1\right)R+\frac{4(\delta b_j+2)}{1-\delta}\frac{c}{H}}{\frac{4(a_j-1)}{1-\delta}R+\frac{4(b_j+2)}{1-\delta}\frac{c}{H}}\right)\!\right)^{\frac{q^{j+2}-q}{q-1}} \mathrm{d}\rho  \label{lower bound estimate J j+1}
\end{align} 
 for $t\geqslant A^{-1}((8a_j+1)R+4(2b_j+1) \frac{c}{H})$, where $\bar{B}$ is defined in \eqref{B bar def}. We denote by $\widetilde{J}_{j+1}(t) $ the $\rho$-integral in the left-hand side of the  last inequality.

We proceed with the lower bound estimate for $\widetilde{J}_{j+1}(t) $ by shrinking the domain of integration. Consequently, for $t\geqslant A^{-1}((8a_j+1)R+4(2b_j+1) \frac{c}{H})$ we obtain
\begin{align}
 \widetilde{J}_{j+1}(t)  & \geqslant  \int_{\delta(A(t)-(8a_j+1)R-4(2b_j+1) \frac{c}{H})}^{A(t)-(8a_j+1)R-4(2b_j+1) \frac{c}{H}}  \frac{\rho^{(n-1)[1 -\frac p2]_+}}{\left(A(t)-\rho-R+\tfrac{2c}{H}\right)} \notag \\
 & \qquad \qquad \qquad  \times \left(\ln \! \left(\frac{A(t)-\rho+\left(\frac{4(\delta a_j-1)}{1-\delta}-1\right)R+\frac{4(\delta b_j+2)}{1-\delta}\frac{c}{H}}{\frac{4(a_j-1)}{1-\delta}R+\frac{4(b_j+2)}{1-\delta}\frac{c}{H}}\right)\!\right)^{\frac{q^{j+2}-q}{q-1}} \mathrm{d}\rho \notag \\
  &  \geqslant \delta^{(n-1)[1-\frac{p}{2}]_+} \left(A(t)-(8a_j+1)R-4(2b_j+1)\tfrac{c}{H}\right)^{(n-1)[1-\frac{p}{2}]_+} \widehat{J}_{j+1}(t), \label{lower bound estimate J tilde j+1}
\end{align} where $ \widehat{J}_{j+1}(t)$ denotes the $\rho$-integral with the remaining factors.

Our goal now is to derive a lower bound estimate for $ \widehat{J}_{j+1}(t)$ in a such way that the power of the logarithmic factor is increased by $1$. According to this purpose, we decrease the argument of the logarithmic term as follows
\begin{align*}
& \frac{A(t)-\rho+\left(\frac{(4(\delta a_j-1)}{1-\delta}-1\right)R+\frac{4(\delta b_j+2)}{1-\delta}\frac{c}{H}}{\frac{4(a_j-1)}{1-\delta}R+\frac{4(b_j+2)}{1-\delta}\frac{c}{H}} \\
& \qquad \qquad \geqslant \frac{A(t)-\rho+\left(\frac{4(\delta a_j-1)}{1-\delta}-1-4a_j\right)R+\left(\frac{4(\delta b_j+2)}{1-\delta}-4(b_j+1)\right)\frac{c}{H}}{\frac{4(a_j-1)}{1-\delta}R+\frac{4(b_j+2)}{1-\delta}\frac{c}{H}} \\
& \qquad \qquad= \frac{A(t)-\rho+\frac{4(2\delta-1) a_j-5+\delta}{1-\delta} R+\frac{4(2\delta-1) b_j+4+4\delta }{1-\delta}\frac{c}{H}}{\frac{4(a_j-1)}{1-\delta}R+\frac{4(b_j+2)}{1-\delta}\frac{c}{H}}
\end{align*} and we use the inequality
\begin{align} \label{upper bound denominator}
A(t)-\rho-R+\tfrac{2c}{H}\leqslant A(t)-\rho+\tfrac{4(2\delta-1) a_j-5+\delta}{1-\delta} R+\tfrac{4(2\delta-1) b_j+4+4\delta }{1-\delta}\tfrac{c}{H}.
\end{align} In particular, we can guarantee the validity of \eqref{upper bound denominator} if we assume  that $\delta> \frac{1}{2}$ and that $a_j,b_j$ fulfill the following conditions
\begin{align}\label{condition aj}
a_j\geqslant \frac{1}{2\delta-1}, \qquad b_j\geqslant - \frac{(1+3\delta)}{2(2\delta-1)}.
\end{align} The condition for $b_j$ in \eqref{condition aj} is trivially true for any $j\in\mathbb{N}$, since $\{b_j\}_{j\in\mathbb{N}}$ is a sequence of positive real numbers. On the other hand, the condition for $a_j$ in \eqref{condition aj} is satisfied for any $j\in\mathbb{N}$ if and only if it is satisfied for $j=0$, being $\{a_j\}_{j\in\mathbb{N}}$ an increasing sequence. Therefore, we impose the following further condition on $a_0$
\begin{align}\label{condition a0 inductive step}
a_0\geqslant \frac{1}{2\delta-1}.
\end{align}
 This condition, together with the one in \eqref{condition a0, b0 first step}, provide the range for $a_0$ that makes our inductive argument successful.
 In conclusion, we proceed with the estimate from below for $\widehat{J}_{j+1}(t)$ as follows
 \begin{align*}
& \widehat{J}_{j+1}(t)  \geqslant  \int_{\delta(A(t)-(8a_j+1)R-4(2b_j+1) \frac{c}{H})}^{A(t)-(8a_j+1)R-4(2b_j+1) \frac{c}{H}}  \left(A(t)-\rho+\tfrac{4(2\delta-1) a_j-5+\delta}{1-\delta} R+\tfrac{4(2\delta-1) b_j+4+4\delta }{1-\delta}\tfrac{c}{H}\right)^{-1}\\
 & \qquad \times \left(\ln \! \left(\frac{A(t)-\rho+\frac{4(2\delta-1) a_j-5+\delta}{1-\delta} R+\frac{4(2\delta-1) b_j+4+4\delta }{1-\delta}\frac{c}{H}}{\frac{4(a_j-1)}{1-\delta}R+\frac{4(b_j+2)}{1-\delta}\frac{c}{H}}\right)\!\right) ^{\frac{q^{j+2}-q}{q-1}} \mathrm{d}\rho  \\
 & \ = \tfrac{q-1}{q^{j+2}-1} \left(\!\ln \! \left(\!\frac{A(t)+\left(\frac{4(2\delta-1) a_j-5+\delta}{(1-\delta)^2}+\frac{\delta(8a_j+1)}{1-\delta}\right) \! R+\left(\frac{4(2\delta-1) b_j+4+4\delta }{(1-\delta)^2}+\frac{4\delta (2b_j+1)}{1-\delta}\right)\!\frac{c}{H}}{\frac{4(a_j-1)}{(1-\delta)^2}R+\frac{4(b_j+2)}{(1-\delta)^2}\frac{c}{H}}\right)\!\!\right)^{\frac{q^{j+2}-1}{q-1}}
 \end{align*}  $t\geqslant A^{-1}((8a_j+1)R+4(2b_j+1) \frac{c}{H})$.
 Having in mind \eqref{||v(t)||^p Lp lower bound step j} for $j+1$, according to the terms appearing in the denominator of the argument of the logarithmic term in the previous estimate, we set 
 \begin{align*}
 \tfrac{a_{j+1}-1}{1-\delta}=\tfrac{4(a_{j}-1)}{(1-\delta)^2}, \qquad \tfrac{b_{j+1}+2}{1-\delta}=\tfrac{4(b_{j}+2)}{(1-\delta)^2}.
 \end{align*} The previous two conditions are equivalent to define
 \begin{align}\label{iterative definition aj, bj}
 a_{j+1}= \tfrac{4}{1-\delta}a_j-\tfrac{3+\delta}{1-\delta}, \qquad  b_{j+1}= \tfrac{4}{1-\delta}b_j+\tfrac{2(3+\delta)}{1-\delta}.
 \end{align} Using iteratively \eqref{iterative definition aj, bj}, we find exactly the representations for $a_j$ and $b_j$ given in \eqref{def aj geometric sequence} and \eqref{def bj geometric sequence}, respectively. 
 
Moreover, by straightforward computations we find that the conditions
\begin{align*}
\begin{cases}
\frac{4(2\delta-1) a_j-5+\delta}{(1-\delta)^2}+\frac{\delta(8a_j+1)}{1-\delta}\geqslant \frac{\delta a_{j+1}-1}{1-\delta} \\
\frac{4(2\delta-1) b_j+4+4\delta }{(1-\delta)^2}+\frac{4\delta (2b_j+1)}{1-\delta} \geqslant \frac{\delta b_{j+1}+2}{1-\delta}
\end{cases}
\end{align*}  are satisfied if and only if the inequalities in \eqref{condition aj} hold. Therefore, thanks to \eqref{condition a0 inductive step}, we can continue the lower bound estimate for $\widehat{J}_{j+1}(t)$, obtaining $t\geqslant A^{-1}((8a_j+1)R+4(2b_j+1) \frac{c}{H})$ 
\begin{align*}
\widehat{J}_{j+1}(t)  & \geqslant  (q-1) q^{-(j+2)} \left(\ln \! \left(\frac{A(t)+\frac{\delta a_{j+1}-1}{1-\delta}  R+\frac{\delta b_{j+1}+2}{1-\delta}\frac{c}{H}}{\frac{a_{j+1}-1}{1-\delta}R+\frac{b_{j+1}+2}{1-\delta}\frac{c}{H}}\right)\!\right)^{\frac{q^{j+2}-1}{q-1}}.
\end{align*}
Combining the lower bound estimates for $\widehat{J}_{j+1}(t)$ and $\widetilde{J}_{j+1}(t)$, we get
\begin{align*}
\widetilde{J}_{j+1}(t) & \geqslant \delta^{(n-1)[1-\frac{p}{2}]_+} (q-1) q^{-(j+2)} \left(A(t)-a_{j+1} R-b_{j+1}\tfrac{c}{H}\right)^{(n-1)[1-\frac{p}{2}]_+} \\ 
& \qquad \qquad  \times  \left(\ln \! \left(\frac{A(t)+\frac{\delta a_{j+1}-1}{1-\delta}  R+\frac{\delta b_{j+1}+2}{1-\delta}\frac{c}{H}}{\frac{a_{j+1}-1}{1-\delta}R+\frac{b_{j+1}+2}{1-\delta}\frac{c}{H}}\right)\!\right)^{\frac{q^{j+2}-1}{q-1}},
\end{align*} where we used the relations $a_{j+1}\geqslant 8 a_j+1$, $b_{j+1}\geqslant 4( 2b_j+1)$ (which turn out to be equivalent to the conditions in \eqref{condition aj}) to lower the first time-dependent factor on the right-hand side that comes from \eqref{lower bound estimate J tilde j+1}.

Finally, putting together the last estimate for $\widetilde{J}_{j+1}(t)$, \eqref{lower bound estimate J j+1} and \eqref{||v(t)||^p Lp lower bound step j intermediate}, we conclude
\begin{align}
\|v(t,\cdot)\|^p_{L^p(\mathbb{R}^n)} 
&  \geqslant     B_{j+1} \varepsilon^{pq^{j+2}}  (1+t)^{\varsigma p \frac{q^{j+2}-1}{q-1}}  \mathrm{e}^{-\frac{1}{2}(b+H)pt} \frac{\left(A(t)-a_{j+1}R-b_{j+1}\tfrac{c}{H}\right)^{(n-1)[1-\frac{p}{2}]_+}}{(A(t)+R)^{(n-1)[1 -\frac p2]_-}} \notag \\
& \qquad \qquad \times \left(\ln \! \left(\frac{A(t)+\frac{\delta a_{j+1}-1}{1-\delta}  R+\frac{\delta b_{j+1}+2}{1-\delta}\frac{c}{H}}{\frac{a_{j+1}-1}{1-\delta}R+\frac{b_{j+1}+2}{1-\delta}\frac{c}{H}}\right)\!\right)^{\frac{q^{j+2}-1}{q-1}}, \label{||v(t)||^p Lp lower bound step j+1}
\end{align} where 
\begin{align} \label{recursive from Bj to B j+1}
B_{j+1}\doteq  D q^{-(j+1)} B_j^q,
\end{align}  with $D= D(n,c,H,b,m^2,\beta,p,a_0,\delta) \doteq K \bar{B} \widehat{C}^p \delta^{(n-1)[1-\frac{p}{2}]_+}  (q-1)q^{-1} $.

 Hence we proved \eqref{||v(t)||^p Lp lower bound step j} for $j+1$ which is exactly \eqref{||v(t)||^p Lp lower bound step j+1} with $B_{j+1}$ given by \eqref{recursive from Bj to B j+1}. Finally, it is convenient for our future considerations to derive an explicit representation of $\ln B_j$. Applying the logarithmic function to both sides of \eqref{recursive from Bj to B j+1} and using the resulting identity in an iterative way, we find
\begin{align}
\ln B_j & = q\ln B_{j-1} -j \ln q+\ln D = q^2\ln B_{j-2} -(j+(j-1)q) \ln q+(1+q)\ln D \notag \\
& = \cdots =  q^j\ln B_{0} -\left(\ \sum_{k=0}^{j-1}(j-k)q^k\right) \ln q+\left(\ \sum_{k=0}^{j-1} q^k\right)\ln D \notag \\
&=  q^j\ln B_{0} -\frac{1}{q-1}\left(\frac{q^{j+1}-1}{q-1}-(j+1)\right) \ln q+\frac{q^j-1}{q-1}\ln D \notag \\
&= q^j \left(\ln B_0 -\frac{q\ln q}{(q-1)^2}+\frac{\ln D}{q-1}\right) +(j+1)\frac{\ln q}{q-1}+\frac{\ln q}{(q-1)^2}-\frac{\ln D}{q-1} \notag \\
& = q^j \ln E +(j+1)\frac{\ln q}{q-1}+\frac{\ln q}{(q-1)^2}-\frac{\ln D}{q-1}, \label{ln Bj representation}
\end{align} where $E\doteq B_0 q^{-q/(q-1)^2} D^{1/(q-1)}$.

\subsection{Improved lower bound estimates for the spatial average of the solution} \label{Subsection sequence lb V(t)}

In the previous subsection we derived the sequence of lower bound estimates for $\|v(t,\cdot)\|^p_{L^p(\mathbb{R}^n)}$ in \eqref{||v(t)||^p Lp lower bound step j}. Now, we are going to use \eqref{||v(t)||^p Lp lower bound step j} to derive a sequence of lower bound estimates for the spatial average $V(t)$.

Indeed, plugging \eqref{||v(t)||^p Lp lower bound step j} in \eqref{fundamental ineq V(t)}, we have
\begin{align*}
V(t)  \geqslant  B_j^{\beta+1} \varepsilon^{q^{j+2}}\mathrm{e}^{-\alpha_2 t} & \int_{A^{-1} (a_jR+b_j\frac{c}{H})}^t \mathrm{e}^{(\alpha_2-\alpha_1) s} \int_{A^{-1} (a_jR+b_j\frac{c}{H})}^s \mathrm{e}^{(\alpha_1-\frac{1}{2}(b+H)q) \tau} \, \Gamma(\tau) \,   \notag  \\ &  \qquad \qquad  \times (1+\tau)^{\varsigma  \frac{q^{j+2}-q}{q-1}}    \frac{(A(\tau)-a_jR-b_j\tfrac{c}{H})^{(n-1)(\beta+1)[1 -\frac p2]_+}}{(A(\tau)+R)^{(n-1)(\beta+1)[1 -\frac p2]_-}} \notag  \\ &  \qquad \qquad \times \left(\ln \! \left(\frac{A(\tau)+\frac{\delta a_j-1}{1-\delta}R+\frac{\delta b_j+2}{1-\delta}\frac{c}{H}}{\frac{a_j-1}{1-\delta}R+\frac{b_j+2}{1-\delta}\frac{c}{H}}\right)\!\right)^{(\beta+1)\frac{q^{j+1}-1}{q-1}}  \! \mathrm{d}\tau  \, \mathrm{d}s
\end{align*} for $t\geqslant A^{-1} (a_jR+b_j\frac{c}{H})$.

Since $A(\tau)-a_jR-b_j\tfrac{c}{H}\geqslant \frac{1}{2}\phi (\tau)$ for $\tau\geqslant A^{-1}(2a_j R+(2b_j+1)\frac{c}{H})$, and $A(\tau)+R\leqslant 2\phi(\tau)$ provided that $\tau \geqslant A^{-1}(R-\frac{2c}{H})$, we have that
\begin{align*}
\frac{(A(\tau)-a_jR-b_j\tfrac{c}{H})^{(n-1)(\beta+1)[1 -\frac p2]_+}}{(A(\tau)+R)^{(n-1)(\beta+1)[1 -\frac p2]_-}} \geqslant 2^{-(n-1)(\beta+1)|1-\frac{p}{2}|}\left(\tfrac{c}{H}\right)^{(n-1)(\beta+1)(1-\frac{p}{2})}\mathrm{e}^{(n-1)H(\beta+1)(1-\frac{p}{2})\tau}
\end{align*} for $\tau\geqslant A^{-1}(2a_j R+(2b_j+1)\frac{c}{H})$. Therefore,
\begin{align*}
V(t) & \geqslant 2^{-(n-1)(\beta+1)|1-\frac{p}{2}|} \mu \left(\tfrac{c}{H}\right)^{(n-1)(\beta+1)(1-\frac{p}{2})} B_j^{\beta+1} \varepsilon^{q^{j+2}}\mathrm{e}^{-\alpha_2 t}  \int_{A^{-1} (2a_jR+(2b_j+1)\frac{c}{H})}^t \mathrm{e}^{(\alpha_2-\alpha_1) s}  \notag  \\ &  \qquad \quad \times \int_{A^{-1} (2a_jR+(2b_j+1)\frac{c}{H})}^s  (1+\tau)^{\varsigma  \frac{q^{j+2}-1}{q-1}} \mathrm{e}^{(\frac{n}{2}+\nu-\frac{1}{p})H\tau}   \notag  \\ &  \quad \qquad \qquad  \times \left(\ln \! \left(\frac{A(\tau)+\frac{\delta a_j-1}{1-\delta}R+\frac{\delta b_j+2}{1-\delta}\frac{c}{H}}{\frac{a_j-1}{1-\delta}R+\frac{b_j+2}{1-\delta}\frac{c}{H}}\right)\!\right)^{(\beta+1)\frac{q^{j+1}-1}{q-1}}   \mathrm{d}\tau  \, \mathrm{d}s
\end{align*} for $t\geqslant A^{-1} (2a_jR+(2b_j+1)\frac{c}{H})$, where we used
\begin{align*}
& \mathrm{e}^{(\alpha_1-\frac{1}{2}(b+H)q+(n-1)(\beta+1)(1-\frac{p}{2})) \tau} \, \Gamma(\tau) \, (1+\tau)^{\varsigma  \frac{q^{j+2}-q}{q-1}} \\ 
& \qquad \qquad =\mu \, (1+\tau)^{\varsigma  \frac{q^{j+2}-1}{q-1}} \mathrm{e}^{(\alpha_1+\varrho_{\mathrm{crit}}-\frac{1}{2}(b+H)q+(n-1)H(\beta+1)(1-\frac{p}{2})) \tau}\\
& \qquad \qquad =\mu \, (1+\tau)^{\varsigma  \frac{q^{j+2}-1}{q-1}} \mathrm{e}^{(\alpha_1+\frac{nH}{2}-\frac{b}{2}-\frac{H}{p})\tau}  =\mu \, (1+\tau)^{\varsigma  \frac{q^{j+2}-1}{q-1}} \mathrm{e}^{(\frac{n}{2}+\nu-\frac{1}{p})H\tau}
\end{align*} with $\varrho_{\mathrm{crit}}$ given by \eqref{def rho crit} and for $\alpha_1,\alpha_2$ we employ the same notations as in Lemma \ref{Lemma ODI}, that is, $\alpha_1=\frac{b}{2}+\nu H$ and $\alpha_2=\frac{b}{2}-\nu H$.
Since we are working with $\varsigma\leqslant 0$, we have
\begin{align}
V(t) & \geqslant  \widetilde{N}  B_j^{\beta+1} \varepsilon^{q^{j+2}}  (1+t)^{\varsigma  \frac{q^{j+2}-1}{q-1}}  \mathrm{e}^{-\alpha_2 t}  \int_{A^{-1} (2a_jR+(2b_j+1)\frac{c}{H})}^t \mathrm{e}^{(\alpha_2-\alpha_1) s}  \notag  \\ &  \quad \times \int_{A^{-1} (2a_jR+(2b_j+1)\frac{c}{H})}^s \mathrm{e}^{(\frac{n}{2}+\nu-\frac{1}{p})H\tau}   \left(\ln \! \left(\frac{A(\tau)+\frac{\delta a_j-1}{1-\delta}R+\frac{\delta b_j+2}{1-\delta}\frac{c}{H}}{\frac{a_j-1}{1-\delta}R+\frac{b_j+2}{1-\delta}\frac{c}{H}}\right)\!\right)^{(\beta+1)\frac{q^{j+1}-1}{q-1}}   \mathrm{d}\tau  \, \mathrm{d}s \label{t integral V(t) lb}
\end{align} for $t\geqslant A^{-1} (2a_jR+(2b_j+1)\frac{c}{H})$, where $\widetilde{N}\doteq 2^{-(n-1)(\beta+1)|1-\frac{p}{2}|} \mu \left(\tfrac{c}{H}\right)^{(n-1)(\beta+1)(1-\frac{p}{2})}$. We focus now on the lower bound estimate for the $\tau$-integral in the right-hand side of the last estimate. For $s\geqslant A^{-1} (4a_jR+2(2b_j+1)\frac{c}{H})$ we may shrink the domain of integration to $[A^{-1}(A(s)/2),s]$, obtaining
\begin{align*}
M(s) & \doteq  \int_{A^{-1}\left(\frac{A(s)}{2}\right)}^s \mathrm{e}^{\left(\frac{n}{2}+\nu-\frac{1}{p}\right)H\tau}   \left(\ln \! \left(\frac{A(\tau)+\frac{\delta a_j-1}{1-\delta}R+\frac{\delta b_j+2}{1-\delta}\frac{c}{H}}{\frac{a_j-1}{1-\delta}R+\frac{b_j+2}{1-\delta}\frac{c}{H}}\right)\!\right)^{(\beta+1)\frac{q^{j+1}-1}{q-1}}   \mathrm{d}\tau \\
&  \geqslant  \left(\ln \! \left(\frac{A(s)+\frac{2(\delta a_j-1)}{1-\delta}R+\frac{2(\delta b_j+2)}{1-\delta}\frac{c}{H}}{\frac{2(a_j-1)}{1-\delta}R+\frac{2(b_j+2)}{1-\delta}\frac{c}{H}}\right)\!\right)^{(\beta+1)\frac{q^{j+1}-1}{q-1}}   \int_{A^{-1}\left(\frac{A(s)}{2}\right)}^s \mathrm{e}^{\left(\frac{n}{2}+\nu-\frac{1}{p}\right)H\tau}   \mathrm{d}\tau  
\\
&  = \frac{\mathrm{e}^{\left(\frac{n}{2}+\nu-\frac{1}{p}\right)H s}}{(\frac{n}{2}+\nu-\frac{1}{p})H}  \left(\ln \! \left(\frac{A(s)+\frac{2(\delta a_j-1)}{1-\delta}R+\frac{2(\delta b_j+2)}{1-\delta}\frac{c}{H}}{\frac{2(a_j-1)}{1-\delta}R+\frac{2(b_j+2)}{1-\delta}\frac{c}{H}}\right)\!\right)^{(\beta+1)\frac{q^{j+1}-1}{q-1}}    \\ &  \qquad \times \left(1-\mathrm{e}^{- \left(\frac{n}{2}+\nu-\frac{1}{p}\right)H\left(s-A^{-1}\left(\frac{A(s)}{2}\right)\right)}\right).
\end{align*} By elementary computations we find,
\begin{align*}
s-A^{-1}\left(\tfrac{A(s)}{2}\right)= -\tfrac 1H \ln \left(1+\mathrm{e}^{-Hs}\right)+ \tfrac 1H  \ln 2,
\end{align*} so that we may estimate for $s\geqslant A^{-1} (4a_jR+2(2b_j+1)\frac{c}{H})$ the last factor in the lower bound for $M(s)$ as follows:
\begin{align*}
1-\mathrm{e}^{- (\frac{n}{2}+\nu-\frac{1}{p})H\left(s-A^{-1}\left(\frac{A(s)}{2}\right)\right)} & = 1-\left(\frac{1+\mathrm{e}^{-H s}}{2}\right)^{\left(\frac{n}{2}+\nu-\frac{1}{p}\right)} \\ & \geqslant 
1-\left(\tfrac{1}{2}+\tfrac{c}{2H} \left(4a_j R+(4b_j+3)\tfrac{c}{H}\right)^{-1}\right)^{\left(\frac{n}{2}+\nu-\frac{1}{p}\right)} \doteq \gamma_j.
\end{align*} Clearly, $\{\gamma_j\}_{j\in\mathbb{N}}$ is an increasing sequence of positive real numbers and $\displaystyle{\lim_{j\to \infty}\gamma_j}=1-2^{-(\frac{n}{2}+\nu-\frac{1}{p})}$.
Combining the previous considerations, we proved that the $s$-integral in \eqref{t integral V(t) lb} can be estimate from below by
\begin{align*}
\widetilde{M}(t) & \doteq  \frac{\gamma_j}{(\frac{n}{2}+\nu-\frac{1}{p})H}   \int_{A^{-1} (4a_jR+2(2b_j+1)\frac{c}{H})}^t \mathrm{e}^{\left(\frac{n}{2}-\nu-\frac{1}{p}\right)H s} \\
& \qquad \qquad \qquad \qquad \qquad \times \left(\ln \! \left(\frac{A(s)+\frac{2(\delta a_j-1)}{1-\delta}R+\frac{2(\delta b_j+2)}{1-\delta}\frac{c}{H}}{\frac{2(a_j-1)}{1-\delta}R+\frac{2(b_j+2)}{1-\delta}\frac{c}{H}}\right)\!\right)^{(\beta+1)\frac{q^{j+1}-1}{q-1}}  \mathrm{d}s,
\end{align*} where we used $\alpha_2-\alpha_1=-2\nu H$.

For $t\geqslant A^{-1} (8a_jR+4(2b_j+1)\frac{c}{H})$, we increase the bottom of the interval of integration so that the domain of the integral is reduced to $[A^{-1}(A(t)/2),t]$. Hence, for $t\geqslant A^{-1} (8a_jR+4(2b_j+1)\frac{c}{H})$ we have
\begin{align*}
\widetilde{M}(t) & \geqslant  \frac{\gamma_j}{(\frac{n}{2}+\nu-\frac{1}{p})H}  \left(\ln \! \left(\frac{A(t)+\frac{4(\delta a_j-1)}{1-\delta}R+\frac{4(\delta b_j+2)}{1-\delta}\frac{c}{H}}{\frac{4(a_j-1)}{1-\delta}R+\frac{4(b_j+2)}{1-\delta}\frac{c}{H}}\right)\!\right)^{(\beta+1)\frac{q^{j+1}-1}{q-1}}   \\
& \qquad \qquad \times \int_{A^{-1}\left(\frac{A(t)}{2}\right)}^t  \mathrm{e}^{\left(\frac{n}{2}-\nu-\frac{1}{p}\right)H s}  \, \mathrm{d}s.
\end{align*} We can estimate the integral in the last inequality for $\widetilde{M}(t)$ analogously to the corresponding term that appeared in the lower bound estimate for $M(s)$, that is, for  $t\geqslant A^{-1} (8a_jR+4(2b_j+1)\frac{c}{H})$
\begin{align*}
\int_{A^{-1}\left(\frac{A(t)}{2}\right)}^t \! \mathrm{e}^{\left(\frac{n}{2}-\nu-\frac{1}{p}\right)H s}  \, \mathrm{d}s \geqslant \frac{\widetilde{\gamma}_j }{\big(\frac{n}{2}-\nu-\frac{1}{p}\big)H} \, \mathrm{e}^{\left(\frac{n}{2}-\nu-\frac{1}{p}\right)H t}, 
\end{align*} where
\begin{align*}
\widetilde{\gamma}_j \doteq 1-\left(\tfrac{1}{2}+\tfrac{c}{2H} \left(8a_j R+(8b_j+5)\tfrac{c}{H}\right)^{-1}\right)^{\left(\frac{n}{2}-\nu-\frac{1}{p}\right)}. 
\end{align*}

Summarizing, for $t\geqslant A^{-1} (8a_jR+4(2b_j+1)\frac{c}{H})$ we proved that 
\begin{align*}
V(t) & \geqslant \widehat{N} \gamma_j \, \widetilde{\gamma}_j \, B_j^{\beta+1} \varepsilon^{q^{j+2}} \mathrm{e}^{(\frac{nH}{2}-\frac{b}{2}-\frac{H}{p})t}  (1+t)^{\varsigma  \frac{q^{j+2}-1}{q-1}} \\
& \qquad \qquad \times \left(\ln \! \left(\frac{A(t)+\frac{4(\delta a_j-1)}{1-\delta}R+\frac{4(\delta b_j+2)}{1-\delta}\frac{c}{H}}{\frac{4(a_j-1)}{1-\delta}R+\frac{4(b_j+2)}{1-\delta}\frac{c}{H}}\right)\!\right)^{(\beta+1)\frac{q^{j+1}-1}{q-1}}, 
\end{align*} where $\widehat{N}\doteq  \widetilde{N} \big( (\frac{n}{2}+\nu-\frac{1}{p}) (\frac{n}{2}-\nu-\frac{1}{p})H^2\big)^{-1} $.

\begin{remark} In the previous lower bound estimate for $V(t)$ the logarithmic term provides an improvement of the lower bound in comparison to what we would obtained if we worked with the same approach as in the proof of \cite[Theorem 1.9]{PalTak22}.
\end{remark}

Next, we require a further assumptions on $a_0$. More precisely, we assume that 
\begin{align} \label{condition a0 final lb V(t)}
a_0\geqslant  \frac{1}{\delta}.
\end{align} As a consequence of \eqref{condition a0 final lb V(t)}, we find that $\delta a_j-1\geqslant 0$ for any $j\in\mathbb{N}$, and, consequently,
\begin{align*}
V(t) & \geqslant \widehat{N} \gamma_j \, \widetilde{\gamma}_j \, B_j^{\beta+1} \varepsilon^{q^{j+2}} \mathrm{e}^{(\frac{nH}{2}-\frac{b}{2}-\frac{H}{p})t}  (1+t)^{\varsigma  \frac{q^{j+2}-1}{q-1}} \left(\ln \! \left(\frac{A(t)}{\frac{4(a_j-1)}{1-\delta}R+\frac{4(b_j+2)}{1-\delta}\frac{c}{H}}\right)\!\right)^{(\beta+1)\frac{q^{j+1}-1}{q-1}}
\end{align*}  for $t\geqslant A^{-1} (8a_jR+4(2b_j+1)\frac{c}{H})$.

For $t\geqslant A^{-1}\left(2^4\big(\tfrac{a_j-1}{1-\delta}R+\tfrac{b_j+2}{1-\delta}\tfrac{c}{H}\big)^2\right)$ it holds the following inequality
\begin{align*}
\ln \! \left(\frac{A(t)}{\frac{4(a_j-1)}{1-\delta}R+\frac{4(b_j+2)}{1-\delta}\frac{c}{H}}\right) \geqslant  \frac{1}{2}\ln A(t).
\end{align*} Furthermore, the condition $t\geqslant \frac{2}{H}\ln (\frac{H}{c}+1)$ implies that $\ln A(t) \geqslant \frac{H t}{2}$, while for $t\geqslant 1$ we may estimate $(1+t)^\varsigma\geqslant 2^\varsigma t^\varsigma$. 

Let us introduce 
\begin{align}\label{def sigma j}
\sigma_j\doteq \max \left\{A^{-1} (8a_jR+4(2b_j+1)\tfrac{c}{H}) ,A^{-1}\left(2^4\big(\tfrac{a_j-1}{1-\delta}R+\tfrac{b_j+2}{1-\delta}\tfrac{c}{H}\big)^2\right),\tfrac{2}{H}\ln (\tfrac{H}{c}+1),1\right\}.
\end{align} Combining the previous estimates,  for $t\geqslant \sigma_j$ we obtained the following lower bound estimate for $V(t)$
\begin{align*}
V(t) & \geqslant N Q^{q^{j+1}} \gamma_j \widetilde{\gamma}_j B_j^{\beta+1}\varepsilon^{q^{j+2}} t^{-\frac{\varsigma+\beta+1}{q-1}}t^{\frac{q^{j+1}}{q-1}(\beta+1+\varsigma q)} \mathrm{e}^{(\frac{nH}{2}-\frac{b}{2}-\frac{H}{p})t}, 
\end{align*} where $N\doteq 2^{-\frac{\varsigma}{q-1}} \widehat{N} \left(\tfrac H4\right)^{-\frac{(\beta+1)}{q-1}}$ and $Q\doteq 2^{\frac{\varsigma q}{q-1}} \left(\tfrac H4\right)^{\frac{\beta+1}{q-1}}$.

Let us denote by $$K_j(t,\varepsilon)\doteq N Q^{q^{j+1}} \gamma_j \widetilde{\gamma}_j B_j^{\beta+1}\varepsilon^{q^{j+2}} t^{-\frac{\varsigma+\beta+1}{q-1}}t^{\frac{q^{j+1}}{q-1}(\beta+1+\varsigma q)}$$ the factor that multiplies the exponential term $\mathrm{e}^{(\frac{nH}{2}-\frac{b}{2}-\frac{H}{p})t} $ in the previous lower bound for $V(t)$. With this notation, we can simply write 
\begin{align} \label{lower bound V in ODI lemma}
V(t)\geqslant K_j(t,\varepsilon) \, \mathrm{e}^{(\frac{nH}{2}-\frac{b}{2}-\frac{H}{p})t} \quad \mbox{for} \ t\geqslant \sigma_j.
\end{align}  Besides, by using the support condition \eqref{support condition sol} and H\"older's inequality, from \eqref{V''+bV'+m2 V} we find that $V$ satisfies
\begin{align}
V''(t)+bV'(t)+m^2V(t) & \gtrsim \Gamma(t) \mathrm{e}^{-nH(\beta+1)(p-1)t} (V(t))^q \notag \\
& \gtrsim \mu (1+t)^\varsigma \mathrm{e}^{\left(\varrho_{\mathrm{crit}}-nH(\beta+1)(p-1)\right)t} (V(t))^q \notag \\
& =  (1+t)^\varsigma \mathrm{e}^{-\left(\frac{nH}{2}-\frac{b}{2}-\frac{H}{p}\right)(q-1)t} (V(t))^q \quad \mbox{for} \ t\geqslant 0, \label{ODI for V in ODI lemma}
\end{align} where we used again \eqref{def rho crit} in the last step.

We point out that Lemma \ref{Lemma ODI} does not cover the case when a polynomial factor appears on right-hand side of the ODI in \eqref{ODI for G statement}. For this reason we consider the next lemma, whose proof is completely analogous to the one of Lemma \ref{Lemma ODI}.

\begin{lemma} \label{Lemma ODI poly} Let $b,m^2$ be nonnegative real numbers such that $b^2\geqslant 4m^2$. We consider the same notations for $\alpha_1$ and $\alpha_2$ as in the statement of Lemma \ref{Lemma ODI}.

 Let $q>1$, $k_0,k_1\in\mathbb{R}$ satisfying 
 \eqref{critical condition k0, k1}, \eqref{condition k0, k1} and  $\ell_0<0$, $\ell_1\in\mathbb{R}$ such that
\begin{align}
& \ell_0+(q-1)\ell_1\geqslant 0. \label{condition l0, l1}
\end{align}
Suppose that $G\in\mathcal{C}^2([0,T))$ satisfies
\begin{align}
& G''(t)+b\, G'(t)+m^2G(t)\geqslant B \, (1+t)^{\ell_0} \mathrm{e}^{k_0 t}|G(t)|^q \qquad \mbox{for} \ t\geqslant 0, \label{ODI for G statement poly} \\
& G(t)\geqslant K (1+t)^{\ell_1}  \mathrm{e}^{k_1 t} \qquad \qquad \qquad \qquad \qquad \qquad \quad \ \, 	\, \mbox{for} \ t\geqslant T_0,  \label{lower bound G statement poly} 
\end{align} and \eqref{Initial condition G statement},  with $T_0\in[0,T)$ and for some  positive constants $B,K$.
Let us define $T_1$ and $K_0$ as in \eqref{def T1 and K0}, where $\vartheta\in (0,\frac{q-1}{2})$ is arbitrarily chosen so that $2\vartheta \ell_1\leqslant \ell_0+(q-1)\ell_1$.

 If $K\geqslant K_0$, then, the lifespan of $G$ is finite and fulfills $T\leqslant 2T_1$.
\end{lemma}

\begin{remark} We emphasize that thanks to the factor $K_j(t,\varepsilon)$ we have a lot of freedom in the choice of the parameter $\ell_1$ in \eqref{lower bound G statement poly}.
\end{remark}

In order to control more accurately $K_j(t,\varepsilon)$, we rewrite it as follows:
\begin{align*}
K_j(t,\varepsilon) & =\exp\left\{ q^{j+1}\left( \ln \varepsilon^q+\ln t^{\frac{\beta+1+\varsigma q}{q-1}}+\ln Q\right) +(\beta+1)\ln B_j+\ln ( \gamma_j \widetilde{\gamma}_j) -\tfrac{\varsigma+\beta+1}{q-1}\ln  t +\ln N\right\}\\
& =\exp\Big\{ q^{j+1}\left( \ln \varepsilon^q+\ln t^{\frac{\beta+1+\varsigma q}{q-1}}+\ln \big(QE^{\frac{1}{p}}\big)\right) +(j+1)\tfrac{(\beta+1)\ln q}{q-1}+\ln ( \gamma_j \widetilde{\gamma}_j) -\tfrac{\varsigma+\beta+1}{q-1}\ln  t \\ & \qquad \qquad +\ln N +\tfrac{(\beta+1)\ln q}{(q-1)^2}-\tfrac{(\beta+1)\ln D}{q-1}\Big\},
\end{align*}
where we used \eqref{ln Bj representation} in the second equality.
Let us introduce the function
\begin{align}\label{def L(t,epsilon)}
L(t,\varepsilon)\doteq \ln\left( \varepsilon^q QE^{\frac{1}{p}}  t^{\frac{\beta+1+\varsigma q}{q-1}}\right).
\end{align} We point out that $L(t,\varepsilon)\geqslant 1$ if and only if $$t\geqslant T_0(\varepsilon)\doteq  E_1 \varepsilon ^{-\frac{p(q-1)}{1+\varsigma p}},$$ where $E_1\doteq \big(\mathrm{e} Q^{-1}E^{-\frac{1}{p}}\big)^{\frac{q-1}{\beta+1+\varsigma q}}$.

Our goal is to apply Lemma \ref{Lemma ODI} or Lemma \ref{Lemma ODI poly} to $V$ (depending on whether $\varsigma=0$ or $\varsigma<0$). We underline that $V$ is twice continuously differentiable (see Equation (3.1) in \cite{PalTak22}). Clearly,  \eqref{ODI for V in ODI lemma} and \eqref{lower bound V in ODI lemma} correspond to the conditions \eqref{ODI for G statement poly} and \eqref{lower bound G statement poly}, respectively, with $\ell_0=\varsigma$. Concerning $\ell_1$, we have some freedom in its choice thanks to the factor $K_j(t,\varepsilon)$. More precisely, we choose a $\ell_1>0$ such that $\varsigma+(q-1)\ell_1>0$.

In particular, for $t\geqslant \max\{T_0(\varepsilon),\sigma_j\}$ it holds
\begin{align*}
K_j(t,\varepsilon) (1+t)^{-\ell_1} & \geqslant \exp\Big\{ q^{j+1} +(\beta+1)(j+1)\tfrac{\ln q}{q-1}+\ln ( \gamma_j \widetilde{\gamma}_j) -\tfrac{\varsigma+\beta+1}{q-1}\ln  t -\ell_1 \ln (1+t) \\ & \qquad \qquad   + \ln N +\tfrac{(\beta+1)\ln q}{(q-1)^2}-\tfrac{(\beta+1)\ln D}{q-1}\Big\}.
\end{align*}
Let us introduce now the family of intervals $\{\mathcal{I}(j)\}_{j\in\mathbb{N}}$, where $\mathcal{I}(j)\doteq [\sigma_j,\sigma_{j+1}]$.

 From \eqref{def sigma j} we see that there exists $j_0=j_0(R,\tfrac{c}{H},a_0,b_0,\delta)\in\mathbb{N}$ such that for $j\geqslant j_0$ we have $\sigma_j=A^{-1}\left(2^4\big(\tfrac{a_j-1}{1-\delta}R+\tfrac{b_j+2}{1-\delta}\tfrac{c}{H}\big)^2\right)$. For $j\geqslant j_0$ and for any $t\in \mathcal{I}(j)$ such that $t\geqslant T_0(\varepsilon)$, we have 
\begin{align}
K_j(t,\varepsilon) (1+t)^{-\ell_1} & \geqslant \exp\Big\{ q^{j+1} +(\beta+1)(j+1)\tfrac{\ln q}{q-1}+\ln ( \gamma_j \widetilde{\gamma}_j) + \ln N +\tfrac{(\beta+1)\ln q}{(q-1)^2}-\tfrac{(\beta+1)\ln D}{q-1}\notag\\ & \qquad \qquad -\left(\tfrac{\varsigma+\beta+1}{q-1}+\ell_1\right)\ln \left(1+ A^{-1}\left(2^4\big(\tfrac{a_{j+1}-1}{1-\delta}R+\tfrac{b_{j+1}+2}{1-\delta}\tfrac{c}{H}\big)^2\right)\right)  \Big\}. \label{lower bound Kj(t,epsilon)}
\end{align} 
We consider now $K_0$ and $\widehat{T}_0(b,m^2,v_0,v_1)$ introduced in Lemmas \ref{Lemma ODI} and \ref{Lemma ODI poly} relatively to $V$. More precisely, we use this lemma with $k_1=\frac{nH}{2}-\frac{b}{2}-\frac{H}{p}$, $k_0=-(q-1)k_1$, $\ell_0=\varsigma$ and $\ell_1$ as discussed above. We stress that the condition \eqref{condition k0, k1} is exactly \eqref{1st lb for V from V0 is dominant} in this case. 

Since the lower bound for $K_j(t,\varepsilon) $ in \eqref{lower bound Kj(t,epsilon)} is divergent as $j\to \infty$, we can fix an index  $J~=~J(n,c,H,R,b,m^2,p,\beta,\sigma,\mu,a_0,b_0,\delta)\in\mathbb{N}$, satisfying $J\geqslant j_0$ in a such way that for any $t\geqslant  A^{-1}\left(2^4\big(\tfrac{a_J-1}{1-\delta}R+\tfrac{b_J+2}{1-\delta}\tfrac{c}{H}\big)^2\right)$ such that $t\geqslant T_0(\varepsilon)$ it holds $K_j(t,\varepsilon)\geqslant K_0$.

We may fix now a $\varepsilon_0=\varepsilon_0(n,c,H,b,m^2,\beta,p,\mu,v_0,v_1,R,a_0,b_0,\delta,J)$ such that 
\begin{align}\label{condition epsilon0}
T_0(\varepsilon_0)& \geqslant \max\left\{A^{-1}\left(2^4\big(\tfrac{a_J-1}{1-\delta}R+\tfrac{b_J+2}{1-\delta}\tfrac{c}{H}\big)^2\right), \widetilde{T}_0,\big(\big(\tfrac{n}{2}-\nu-\tfrac{1}{p}\big)H\big)^{-1} \right\},
\end{align} due to the fact that all terms on the right-hand side are independent of $\varepsilon$.
Consequently, for any $\varepsilon\in(0,\varepsilon_0]$, since $T_0(\varepsilon)\geqslant T_0(\varepsilon_0)$, we have $V(t)\geqslant K_0 (1+t)^{\ell_1} \mathrm{e}^{(\frac{nH}{2}-\frac{b}{2}-\frac{H}{p})t} $ for any $t\geqslant T_0(\varepsilon)$.
In conclusion, Lemmas \ref{Lemma ODI} and \ref{Lemma ODI poly} provide the following upper bound for the lifespan $T(\varepsilon)$ of $V$
\begin{align*}
T(\varepsilon) & \leqslant 2\max\left\{T_0(\varepsilon),\widetilde{T}_0,\big(\big(\tfrac{n}{2}-\nu-\tfrac{1}{p}\big)H\big)^{-1} \right\}\leqslant 2T_0(\varepsilon) = 2E_1 \varepsilon ^{-\frac{p(q-1)}{1+\varsigma p}},
\end{align*} where in the second inequality we used \eqref{condition epsilon0}. This concludes the proof of Theorem \ref{Thm anti dS nlin lb poly growth} for $\varsigma\in (-\frac{1}{p},0]$.

\subsection{Upper bound estimates for the lifespan when $\varsigma >0$} \label{Subsubsection varsigma >0}

In the previous proof we showed the validity of \eqref{upper bound for the lifespan polynomial growth anti dS nlin} when $\varsigma\in (-\frac{1}{p},0]$. Of course, we can repeat the same argument as before when $\varsigma>0$ obtaining the same upper bound estimate for the lifespan as in the case $\varsigma=0$, that is $T(\varepsilon)\lesssim \varepsilon ^{-p(q-1)}$.  On the other hand, when $\varsigma>0$ and $\varrho=\varrho_{\mathrm{crit}}$ we may use the same approach as in the proof of \cite[Theorem 1.9]{PalTak22} to prove the blow-up in finite time of $V$ and the lifespan estimate $T(\varepsilon)\lesssim \varepsilon^{-\frac{q-1}{\varsigma}}$. Comparing the last two lifespan estimates when $\varsigma>0$, we conclude the validity of \eqref{upper bound for the lifespan polynomial growth anti dS nlin}.

\section{Final remarks and open problems}

We point out explicitly that we expect that the lifespan estimate in \eqref{upper bound for the lifespan polynomial growth anti dS nlin} is not sharp when $\varsigma>0$.

 We stress also that when $n>N$ in the double limit case $\varrho=\varrho_{\mathrm{crit}}(n,H,b,m^2,\beta,p)$ and $\varsigma=\varsigma_{\mathrm{crit}}(n,H,b,m^2,p)$ we are not able to prove the blow-up in finite time of $V$. This is due to the fact that the function $L(t,\varepsilon)$ does no longer depend on $t$ as the power for $t$ is $0$ in \eqref{def L(t,epsilon)} when $\varsigma=\varsigma_{\mathrm{crit}}(n,H,b,m^2,p)$.


\section*{Acknowledgments}

A. Palmieri is supported by the\emph{ Japan Society for the Promotion of Science} (JSPS) – JSPS Postdoctoral Fellowship for Research in Japan (Short-term) (PE20003) – and is member of the \emph{Gruppo Nazionale per L’Analisi Matematica, la Probabilità e le loro Applicazioni} (GNAMPA) of the \emph{Instituto Nazionale di Alta Matematica} (INdAM). H. Takamura is partially supported by the Grant-in-Aid for Scientific Research (B) (No.18H01132), \emph{Japan Society for the Promotion of Science}.

\end{document}